\documentclass[12pt]{amsart}

\newtheorem{thm}{Theorem}[section]

\newtheorem{lem}[thm]{Lemma}

\newtheorem{prp}[thm]{Proposition}

\begin{document}
\title{On the density of primes with a set of quadratic residues or non-residues in given arithmetic progression}
\author{Steve 
wright} 
\address{Steve Wright, Department of Mathematics and Statistics,
Oakland University, Rochester, MI 48309. Email: wright@oakland.edu}

\begin{abstract}
Let $\mathcal{A}$ denote a finite set of arithmetic progressions of positive integers and let $s \geq 2$ be an integer. If the cardinality of $\mathcal{A}$ is at least 2 and $U$ is the union formed from certain arithmetic progressions of length $s$ taken from each element of $\mathcal{A}$, we calculate the asymptotic density of the set of all prime numbers $p$ such that $U$ is a set of quadratic residues of $p$ and the asymptotic density of the set of all primes $p$ such that $U$ is a set of quadratic non-residues of $p$.

 \end{abstract}
 
 \maketitle
\markboth{}{}
\noindent \textit{keywords}: \textrm{quadratic residue, quadratic non-residue, arithmetic progression, density of a set of primes, asymptotic approximation}

\section{Introduction}
\label{intro}
If $p$ is an odd prime, an integer $z$ is said to be a \emph{quadratic residue} (respectively, \emph{quadratic non-residue}) of $p$ if the equation $x^2 \equiv z $ mod $ p$ has (respectively, does not have) a solution $x$ in integers. It is a theorem going all the way back to Euler that exactly half of the integers from 1 through $p-1$ are quadratic residues of $p$, and it is a fascinating problem to investigate the various ways in which these residues are distributed among 1, 2,\dots, $p-1$. In this paper, our particular interest lies in measuring the size of the set of odd primes $p$ such that $p$ has a set of quadratic residues or non-residues that form a union of two or more given arithmetic progressions all of a fixed length.

We begin with a litany of notation and terminology that will be used systematically throughout the rest of this paper. If $m \leq n$ are integers, then $[m,n]$ will denote the set of all integers that are at least $m$ and no greater than $n$, listed in increasing order, and $[m,+\infty)$ will denote the set of all integers that exceed $m-1$, also listed in increasing order. If $\{a(p)\}$ and $\{b(p)\}$ are sequences of real numbers defined for all primes $p$ in an infinite set $S$, then we will say that $a(p)$ is (sharply) \emph{asymptotic to} $b(p)$ \emph{as}  $p \rightarrow +\infty$ \emph{inside} $S$, denoted as $a(p) \sim b(p)$, if
\[
\lim_{\substack{p\rightarrow +\infty \\ p\in S}} \frac{a(p)}{b(p)} = 1,
\]
and if $S = [1,+\infty)$, we simply delete the phrase ``inside $S$". If $A$ is a set then $|A|$ will denote the cardinality of $A$, $ 2^A$ will denote the set of all subsets of $A$, $\mathcal{E}(A)$ will denote the set of all nonempty finite subsets of $A$ of even cardinality, and $\emptyset$ will denote the empty set. Finally, we note once and for all that $p$ will always denote a generic odd prime.

In order to state precisely what we wish to do here, it will be convenient to recall one of the principal results from [2]. Let  $(m, s)\in [2, +\infty)\times [1, +\infty)$, let $\textbf{a}=(a_1,\dots, a_m)$, (respectively, $\textbf{b}=(b_1,\dots, b_m)$) be an $m$-tuple of nonnegative (respectively, positive) integers such that $(a_i, b_i)\not= (a_j, b_j)$ for $i\not= j$, and let  $(\textbf{a}, \textbf{b})$ denote the $2m$-tuple $(a_1,\dots, a_m, b_1,\dots, b_m)$. We will refer to $(\textbf{a}, \textbf{b})$ as a \emph{standard} $2m$-tuple.We then let $AP(\textbf{a},\textbf{b};s)$ denote the set
\[
\Big\{\bigcup_{j=1}^{m} \{a_j+b_j(n+i):i\in [0,s-1]\}:n\in [1,+\infty)\Big\}.
\]
If $AP(a_i, b_i)$ denotes the arithmetic progression
\[
\{a_i+b_in: n \in [1, +\infty)\}
\]
with initial term $a_i$ and difference $b_i$, $i \in [1, m]$, then the elements of $AP(\textbf{a},\textbf{b};s)$ consist precisely of the sets formed by taking an $n \in [1, +\infty)$, choosing from $AP(a_i, b_i)$ the arithmetic progression with initial term $a_i+b_in$ and length $s$, $i \in [1, m]$, and then taking the union of  all of these progressions of length $s$. We are interested in the primes $p$ such that an element of $AP(\textbf{a},\textbf{b};s)$ is either a set of quadratic residues of $p$ or, respectively, a set of quadratic non-residues of $p$.

If $p$ is an odd prime and $\mathbb{Z}_p$ is the field of $p$ elements, then the Legendre symbol of $p$ defines a real (primitive) multiplicative character $\chi_p : \mathbb{Z}_p \rightarrow [-1,1]$ on $\mathbb{Z}_p$.  We take $\varepsilon \in \{-1, 1\}$, let
\begin{equation*}
q_\varepsilon(p)=|\{A\in AP(\textbf{a}, \textbf{b}; s)\cap 2^{[1, p-1]}: \chi_p(a)=\varepsilon, \ \textrm{for all}\ a\in A\}|,
\end{equation*} 
and note that the value of $q_\varepsilon(p)$ for $\varepsilon=1$ (respectively, $\varepsilon=-1$) counts the number of elements of $AP(\textbf{a}, \textbf{b}; s)$ that are sets of quadratic residues (respectively, non-residues) of $p$ that are located inside $[1, p-1]$. 

In [2], the sharp asymptotic behavior of  $q_{\varepsilon}(p)$ as $p \rightarrow  +\infty$ was determined. It transpires that
$q_{\varepsilon}(p)$ either has an asymptotic limit as $p \rightarrow  +\infty$  or $q_{\varepsilon}(p)$ asymptotically oscillates infinitely often between $0$ and an asymptotic limit as $p \rightarrow  +\infty$ through a certain infinite set of primes. In order to more precisely describe this behavior, several ingredients from an appropriate recipe must first be listed. Begin by considering the set $B$ of \emph{distinct} values of the coordinates of $\textbf{b}$. If we declare the coordinate $a_i$ of $\textbf{a}$ and the coordinate $b_i$ of $\textbf{b}$ to \emph{correspond} to each other, then for each $b\in B$, we let $A(b)$ denote the set of all coordinates of $\textbf{a}$ whose corresponding coordinate of $\textbf{b}$ is $b$. We then relabel the elements of $B$ as $b_1,\dots,b_k$, say, and for each $i\in [1, k]$, set  
\begin{equation*}
 S_i=\bigcup_{a\in A(b_i)} \{ab_i^{-1}+j: j\in [0, s-1]\}. \tag{1.1} 
 \]
 
The next ingredient is a certain collection of subsets of $[1, k]$ which is constructed from the sets $S_1, \dots, S_k$ in the following manner: let 
\[
\mathcal{K}=\Big\{\emptyset\not= K\subseteq [1, k]: \bigcap_{i\in K}\ S_i\not= \emptyset\Big\},
\]
\[
T(K)=\Big(\bigcap_{i\in K}\ S_i \Big)\cap \Big(\bigcap_{i\in [1, k]\setminus K}\ ( \textbf{Q}\setminus S_i)\Big), K\in \mathcal{K},
\]
where $\textbf{Q}$ denotes the set of all rational numbers, and let
\[
\mathcal{K}_{\max}=\{K\in \mathcal{K}: T(K)\not= \emptyset\}.
\]
The set of subsets of $[1, k]$ that we need is then defined to be the set
\[
\Lambda(\mathcal{K})=\bigcup_{K\in \mathcal{K}_{\max}} \mathcal{E}(K).
\]
N.B. $\Lambda(\mathcal{K})$ is empty if and only if the sets $S_1, \dots, S_k$ are pairwise disjoint. 

Suppose that $\Lambda(\mathcal{K})$ is not empty. It will be convenient to declare that $p$ is an \emph{allowable prime} if no element of $B$ has $p$ as a factor. If $p$ is an allowable prime then we define the $(\textbf{a}, \textbf{b})$-\emph{signature of p} to be the multiset of $\pm 1$'s given by
\[
\Big\{\chi_p\Big(\prod_{i\in I}\  b_i\Big): I\in \Lambda(\mathcal{K}))\Big\}
\]
and then set $\Pi_+(\textbf{a}, \textbf{b})$ (respectively, $ \Pi_-(\textbf{a}, \textbf{b}$)) equal to the set of all allowable primes $p$ such that the $(\textbf{a}, \textbf{b})$-signature of $p$ contains only $1$'s (respectively, contains a $-1$).  For the final ingredients of our recipe, we take
\[
b=\max \{b_1, \dots, b_k\},
\]
\[
\kappa=\Big| \bigcup_{i=1}^k\  S_i \Big|.
\]

The asymptotic behavior of $q_{\varepsilon}(p)$ can now be precisely described. According to Theorem 6.1 of [2],
\vspace{0.3cm}

\noindent (i) if either $S_1, \dots, S_k$ are pairwise disjoint or for all $I\in \Lambda(\mathcal{K}),\  \prod_{i\in I} b_i$ is a square, then
\[
q_\varepsilon(p)\sim (b\cdot 2^\kappa)^{-1}p\ \textrm{as}\ p\rightarrow +\infty, \textrm{or}
\]
\vspace{0.1cm}

\noindent (ii) if there exists $I\in \Lambda(\mathcal{K})$ such that $\prod_{i\in I} b_i$ is not a square, then

\noindent (a)  $\Pi_+(\textbf{a}, \textbf{b})$ and $\Pi_-(\textbf{a}, \textbf{b})$ are both infinite,

\noindent (b) $q_\varepsilon(p)=0$ for all $p$ in $\Pi_-(\textbf{a}, \textbf{b})$, and

\noindent (c) as $p\rightarrow +\infty$ inside $\Pi_+$,
\[
q_\varepsilon(p)\sim (b\cdot2^{\kappa})^{-1}p\ .
\vspace{0.2cm}
\]

The problem that is of interest to us in this article stems from the situation present in statement (ii). In that case, sets of quadratic residues and non-residues form inside $AP(\textbf{a}, \textbf{b}; s)\cap 2^{[1, p-1]}$ only for all primes that are sufficiently large inside $\Pi_+(\textbf{a}, \textbf{b})$, and for no other allowable primes. A natural and interesting question which therefore arises asks: how large can the set $\Pi_+(\textbf{a}, \textbf{b})$ be and how can we measure its size in an accurate way?

A good way to measure the size of an infinite set $\Pi$ of primes is to calculate its (natural or asymptotic) density. If we let P denote the set of all prime numbers then the \emph{density of} $\Pi$ (in P) is defined to be the limit
\[
\lim_{x\rightarrow +\infty } \frac{\big| \{p \in \Pi:\ p \leq x \}\big|}{\big| \{p \in \textrm{P}:\ p \leq x \}\big|} ,
\]
provided that this limit exists. Roughly speaking, the density of $\Pi$ measures the ``proportion" of the set of all primes that are contained in $\Pi$. We will answer the question posed at the end of the previous paragraph by calculating the density of   
 $\Pi_+(\textbf{a}, \textbf{b})$. Because $\Pi_+(\textbf{a}, \textbf{b})$ and $\Pi_-(\textbf{a}, \textbf{b})$ are disjoint sets with only finitely many primes outside of their union, it follows that the density of  $\Pi_-(\textbf{a}, \textbf{b})$ is $1$ minus the density of $\Pi_+(\textbf{a}, \textbf{b})$. Hence a determination of the density of  $\Pi_+(\textbf{a}, \textbf{b})$ also yields a precise measure of the size of the set of primes $p$ such that \emph{no} element of $AP(\textbf{a}, \textbf{b}; s) \cap 2^{[1, p-1]}$ is either a set of quadratic residues or a set of quadratic non-residues of $p$.

As we will see in Lemma 3.1 in section 3, the sets $\{b_i: i \in I\}, I \in \Lambda(\mathcal{K})$, determine the primes in $\Pi_+(\textbf{a}, \textbf{b})$, and so it behooves us to investigate the structure of the set $\Lambda(\mathcal{K})$ in some depth. We do this by deriving in section 2 a very useful combinatorial formula for $\Lambda(\mathcal{K})$, and this formula will hence play an important role in our computation of the density of $\Pi_+(\textbf{a}, \textbf{b})$.  We also present in section 2  some additional combinatorial results and some results on computation of the density that will be required for the calculations performed in section 3. In the latter section, the density of $\Pi_+(\textbf{a}, \textbf{b})$ will be computed for  standard $2m$-tuples $(\textbf{a}, \textbf{b})$ for which certain conditions on the square-free parts of the coordinates of $\textbf{b}$ are satisfied; this is the content of Theorem 3.7, the primary result of this paper.  We proceed by following  a simple strategy: we first decompose $\Pi_+(\textbf{a}, \textbf{b})$ into a finite, pairwise disjoint union of certain sets, then use the density  results of section 2 to calculate the density of these sets, sum everything up, and, finally, use the combinatorial results of section 2 to evaluate this sum. In order to illustrate perspicuously   how Theorem 3.7 determines the density of $\Pi_+(\textbf{a}, \textbf{b})$, we will apply it to the so-called admissible $2k$-tuples. A  standard $2k$-tuple $(\textbf{a}, \textbf{b})$ is  \emph{admissible} if $k \geq2$, the coordinates of $\textbf{b}$ are distinct, and $a_ib_j-a_jb_i \not= 0$ for all $i \not= j$. When $(\textbf{a}, \textbf{b})$ is admissible, we will see that the parameters in the formulae for the density of $\Pi_+(\textbf{a}, \textbf{b})$ that are given by Theorem 3.7 can be calculated by an elegant geometric procedure.

\section{Preliminaries}
In this section we set up the mathematical technology that is required to carry out the calculation of the density of $\Pi_+(\textbf{a}, \textbf{b})$ to be performed in section 3. We begin with two lemmas that will be used to determine the densities of various sets. The first lemma is due to Filaseta and Richman [1, Theorem 2] and the second can be found in [3, Theorem 3.3 ].
\begin{lem}
\label{lem1}
If $S$ is a nonempty finite set of primes and $\varepsilon: S \rightarrow \{-1, 1\}$ is a given function of $S$ into $\{-1, 1\}$ then $2^{-|S|}$ is the density of the set $\{p: \chi_p(z)= \varepsilon(z), \textrm{for all}\  z \in S\}$. 
\end{lem}

The statement of the second lemma requires some preparatory notation. Let $F$ denote the  Galois field $[1, +\infty) / 2[1, +\infty)$ of 2 elements , let $A$ be a finite nonempty subset of $[1, +\infty)$, let $n=|A|$, and let $F^n$ denote the vector space over $F$ of dimension $n$. We arrange the elements $a_1 < \dots < a_n$ of $A$ in increasing order and then define the map $v: 2^{A} \rightarrow F^n$ as follows: if $S \subseteq A$ then the $i$-th coordinate of $v(S)$ is 1 (respectively, 0) if $a_i \in S$ (respectively, $a_i \notin S$). If $z \in [1, +\infty)$, then we denote by $\pi_{\textrm{odd}}(z)$ the set of prime factors of $z$ of odd multiplicity.
\begin{lem}
\label{lem2}
If $S$ is a nonempty finite subset of $[1, +\infty), T=S \setminus\{m^2: m \in [1, +\infty)\}, \mathcal{T}=\{\pi_{\textnormal{odd}}(z): z \in T\}, A= \bigcup \{T: T \in \mathcal{T}\}, n=|A|,$ and 
\begin{center}
$d =$ the dimension of the linear span of $v(\mathcal{T})$ in $F^n$,
\end{center}
then $2^{-d}$ is the density of the set $\{p: \chi_p(z)=1,for\  all\  z \in S\}$.

\end{lem}

The next lemma records some simple enumerative combinatorics that will prove useful in section 3. 
\begin{lem}
\label{lem3}
$([3, Lemma\ 3.2])$ If $A$ is a nonempty finite subset of $[1, +\infty)$, $n=|A|$, $ \mathcal{S}$ and $\mathcal{T}$ are disjoint subsets of $2^A$ and $d$ is the dimension of the linear span of $v(\mathcal{S} \cup \mathcal{T})$ in $F^n$ then the cardinality of the set
\begin{center}
$\{N \subseteq A: |N \cap S|\ is\ odd, for\ all\ S \in \mathcal{S}\ and\ |N \cap T|\ is\ even, for\ all\ T \in \mathcal{T}\}$
\end{center}
is either $0$ or $2^{n-d}$.
\end{lem}

The calculation of the density in section 3 will require a criterion for when the set in the conclusion of Lemma 2.3 is nonempty. In order to state it we recall that the \emph{symmetric difference} $A \triangle B$ of sets $A$ and $B$ is defined as $(A \setminus B) \cup (B \setminus A)$. The symmetric difference operation is commutative and associative, hence if $\{A_1, \dots, A_m\}$ is a finite set of sets then the repeated symmetric difference
\[
A_1 \triangle \cdots \triangle A_m
\]
is unambiguously defined. In fact, one can prove that 
\begin{equation*}
A_1 \triangle \cdots \triangle A_m=\Big\{a \in \bigcup_{i=1}^m\ A_i: |\{A_j: a \in A_j|\ \textrm{is odd} \Big\}. \tag{2.1}
\]
The next lemma is a straightforward reformulation of [3, Proposition 3.5].
\begin{lem}
\label{lem4} If $A$ is a nonempty finite set and $ \mathcal{S}$ and $\mathcal{T}$ are disjoint subsets of $2^A$, with $\emptyset \notin \mathcal{S}$, then the set 
\begin{center}
$\{N \subseteq A: |N \cap S|\ is\ odd, for\ all\ S \in \mathcal{S}\ and\ |N \cap T|\ is\ even, for\ all\ T \in \mathcal{T}\}$
\end{center}
is not empty if and only if for each subset $U$ of $ \mathcal{S} \cup \mathcal{T} \cup \{\emptyset\} $ of odd cardinality, either the cardinality of $U \cap (\mathcal{T} \cup  \{\emptyset\})$ is odd or the repeated symmetric difference of the elements of $U$ is not empty.
\end{lem}

We will now present for the set $\Lambda(\mathcal{K})$ that was defined for standard $2m$-tuples $(\textbf{a}, \textbf{b})$ in section 1 a very useful combinatorial formula. The formula requires the idea of an overlap diagram, defined and studied in [2], and so we will discuss that first.

Begin by choosing $(n, s)\in [1,+\infty)\times [2,+\infty)$ and let $ \textbf{g}=(g(1),\dots,g(n))$ be an $n$-tuple of positive integers. We use \textbf{g} to construct the following array of points. In the plane, place $s$ points horizontally one unit apart, and label the $j$-th point as $(1, j-1)$ for each $j\in [1, s]$. This is \emph{row $1$}. Suppose that row $i$ has been defined. One unit vertically down and $g(i)$ units horizontally to the right of the first point in row $i$, place $s$ points horizontally one unit apart, and label the $j$-th point as $(i+1, j-1)$ for each $j\in [1, s]$. This is $\emph{row}\ i+1$.
The array of points so formed by these $n+1$ rows is called the \emph{overlap diagram of} \textbf{g}, the sequence \textbf{g} is called the \emph{gap sequence} of the overlap diagram, and a nonempty set that is formed by the intersection of the diagram with a vertical line is called a \emph{column} of the diagram. N.B. We do not distinguish between the different possible positions in the plane which the overlap diagram may occupy. A typical example with $n=3, s=8$, and gap sequence (3, 2, 2) looks like
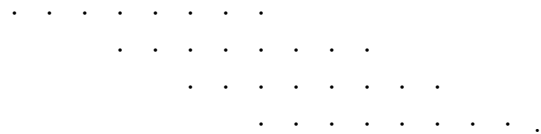
\vspace{0.5cm}
\begin{center}
\begin{tabular}{ccccccccccccccc}
 $\cdot$&$\cdot$&$\cdot$&$\cdot$&$\cdot$&$\cdot$&$\cdot$&$\cdot$&&&&&&&\\

&&&$\cdot$&$\cdot$&$\cdot$&$\cdot$&$\cdot$&$\cdot$&$\cdot$&$\cdot$&&\\
&&&&&$\cdot$&$\cdot$&$\cdot$&$\cdot$&$\cdot$&$\cdot$&$\cdot$&$\cdot$\\
&&&&&&&$\cdot$&$\cdot$&$\cdot$&$\cdot$&$\cdot$&$\cdot$&$\cdot$&$\cdot$\     \ .\\
\end{tabular}
\end{center}
\vspace{0.5cm} 

Next, we need to describe how and where rows overlap in an overlap diagram. Begin by first noticing that if $(g(1),\dots,g(n))$ is the gap sequence, then row $i$ overlaps row $j$ for $i<j$ if and only if
 \[
 \sum_{r=i}^{j-1} g(r)\leq s-1;
 \]                                                         
in particular, row $i$ overlaps row $i+1$ if and only if $g(i)\leq s-1$. Now let $\mathcal{G}$ denote the set of all subsets $G$ of $[1, n]$ such that $G$ is a nonempty set of consecutive integers maximal with respect to the property that $g(i)\leq s-1$ for all $i\in G$. If $\mathcal{G}$ is empty then $g(i)\geq s$ for all $i\in [1, n]$, and so there is no overlap of rows in the diagram. Otherwise there exists $m\in[1, 1+[(n-1)/2]]$ and strictly increasing sequences $(l_1,\dots,l_m)$ and $(M_1,\dots,M_m)$ of positive integers, uniquely determined by the gap sequence of the diagram, such that $l_i\leq M_i$ for all $i\in [1, m], 1+M_i\leq l_{i+1}$ if $i\in [1, m-1]$, and
\[
\mathcal{G}=\{[l_i, M_i]: i\in [1, m]\}.
\]
In fact, $l_{i+1}> 1+M_i$ if $i\in [1, m-1]$, lest the maximality of the elements of $\mathcal{G}$ be violated. It follows that the intervals of integers $[l_i, 1+M_i], i\in [1, m]$, are pairwise disjoint. 

The set $\mathcal{G}$ can now be used to locate the overlap between rows in the overlap diagram like so: for $i\in [1, m]$, let
\[
B_i=[l_i, 1+M_i],
\]
and set
\[
\mathcal{B}_i=\textrm{the set of all points in the overlap diagram whose labels are in}\ B_i\times [0, s-1].
\]
We refer to $\mathcal{B}_i$ as the \emph{i-th block} of the overlap diagram. Thus the blocks of the diagram are precisely the regions in the diagram in which rows overlap. 

Our intent now is to use certain overlap diagrams defined by means of a standard $2m$-tuple to calculate $\Lambda(\mathcal{K})$. We fix a standard $2m$-tuple $(\textbf{a}, \textbf{b})$ and proceed to construct this series of overlap diagrams.

Recall from the  introduction that $B=\{b_1,\dots,b_k\}$ denotes the set of distinct values of the coordinates of $\textbf{b}$, $A(b_i)$ denotes the set of coordinates of $\textbf{a}$ corresponding to $b_i$,
\[
S_i=\bigcup_{a\in A(b_i)} \{ab_i^{-1}+j: j\in [0, s-1]\}, i \in [1, k], 
 \]
\[
\mathcal{K}=\Big\{\emptyset\not= K\subseteq [1, k]: \bigcap_{i\in K}\ S_i\not= \emptyset\Big\},
\]
\[
T(K)=\Big(\bigcap_{i\in K}\ S_i \Big)\cap \Big(\bigcap_{i\in [1, k]\setminus K}\ ( \textbf{Q}\setminus S_i)\Big), K\in \mathcal{K},\ \textrm{and}
\]
\[
\mathcal{K}_{\max}=\{K\in \mathcal{K}: T(K)\not= \emptyset\}.
\]
Let $Q_i=\{a/b_i: a \in A(b_i)\}, i \in [1, k]$, set $Q=\bigcup_i Q_i$ and define the equivalence relation $\approx$ on $Q$ as follows: if $q$ and $q^{\prime}$ are elements of $Q$ then  $q \approx q^{\prime}$ if $q-q^{\prime}$ is an integer. For each equivalence class $E$ of $\approx$, we form the nonempty and pairwise disjoint set of all subsets $R$ of $E$ such that $R$ is maximal with respect to the property that the distance between consecutive elements of $R$ does not exceed $s-1$. We denote by $\mathcal{R}$ the set of all subsets $R$ which arise from all equivalence classes $E$ of $\approx$ via this construction, and observe that the elements of $\mathcal{R}$ form a partition of $Q$.

Let $R \in \mathcal{R}$. We will use $R$ to construct an overlap diagram $\mathcal{D}(R)$. If $|R|=1$, then $\mathcal{D}(R)$ consists of a single row of $s$ points in the plane, spaced one unit apart. Suppose that $\nu=|R| \geq 2$. Arrange the points $r_1< \dots <r_{\nu}$ of $R$ in increasing order, let $d_i=r_{i+1}-r_i$ for $i \in [1, \nu-1]$ and let  $\mathcal{D}(R)$ denote the overlap diagram with gap sequence $(d_1,\dots, d_{\nu})$. Because $d_i \leq s-1$ for all $i \in [1, \nu-1]$, $\mathcal{D}(R)$ consists of a single block.

The family of overlap diagrams $\mathcal{D}(R), R \in \mathcal{R},$ will be used to give the promised formula for $\Lambda(\mathcal{K})$. This will require a certain labeling of the points of these diagrams, which we describe first. N.B. This labeling will in general be different from the labeling of the points of an overlap diagram that was used to define the blocks of the diagram. Let $R \in \mathcal{R}$. The diagram  $\mathcal{D}(R)$ has $|R|$ rows, with each row containing $s$ points. With the elements of $R$ listed in increasing order, for each $i \in [1, |R|]$ we let $r_i$ denote the $i$-th element of $R$. Then proceeding from left to right in row $i$ of $\mathcal{D}(R)$, we take $l \in [1, s]$ and label the $l$-th point of that row as $(r_i, l-1)$. We will say that $r_i$ \emph{labels} row $i$ and also that an element $q$ of $Q$ \emph{corresponds to} row $i$ in $\mathcal{D}(R)$ if $q=r_i$.

Let $C$ be a column of $\mathcal{D}(R)$. We identify $C$ with the subset of $Q \times [0, s-1]$ defined by
\[
\{(q, i)\in Q\times [0, s-1]: (q, i)\ \textrm{is the label of a point in}\ C\},
\]
and then let $\mathcal{C}(R)$ denote the set of all subsets of $Q \times [0, s-1]$  which arise from all such identifications. By a slight abuse of terminology, we will refer to the elements of $\mathcal{C}(R)$  as the columns of  $\mathcal{D}(R)$.  Let
\begin{equation*}
\mathcal{C}= \bigcup_{R \in \mathcal{R}}\ \mathcal{C}(R). \tag{2.2}
\]
If $\theta:Q\times [0, s-1] \rightarrow Q$ is the canonical projection of $Q\times [0, s-1]$ onto $Q$, then for each $C \in \mathcal{C},$ we set
\[
K(C)=\{i \in [1, k]: Q_i \cap \theta(C) \not= \emptyset\}.
\]

Suppose now that the set
\[
\Big\{R \in \mathcal{R}: \mathcal{D}(R)\ \textrm{has a column}\ C \textrm{ such that } |K(C)| \geq 2 \Big\}
\]
is not empty and let $\{R_1,\dots, R_v\}$ be an enumeration of this set. We define the \emph{quotient diagram of} $(\textbf{a}, \textbf{b})$ to be the $v$-tuple of overlap diagrams $(\mathcal{D}(R_1), \dots, \mathcal{D}(R_v))$. The overlap diagrams $\mathcal{D}(R_1), \dots, \mathcal{D}(R_v)$ are called the \emph{blocks} of the quotient diagram.

The following lemma gives a formula which calculates $\Lambda(\mathcal{K})$ by means of the sets $K(C), C \in \mathcal{C}$. We use the notation just introduced in its statement.
\begin{lem}
\label{lem5}
\[
\Lambda(\mathcal{K})= \bigcup_{C \in \mathcal{C}}\ \mathcal{E}(K(C)).
\]
\end{lem}

\emph{Proof}. In order to verify this lemma, it suffices to show that $K \in \mathcal{K}_{\max}$ if and only if there exists $R \in \mathcal{R}$ and a column $C$ of $\mathcal{D}(R)$ such that $K=K(C)$.

For the purposes of the following argument, we fix a horizontal coordinate axis in the plane and locate each overlap diagram $\mathcal{D}(R)$  in the plane so that if $q$ is the smallest element of $R$ then the first point in the first row of $\mathcal{D}(R)$ lies directly over the point $q$ on the coordinate axis.

Let $K \in \mathcal{K}_{\max}$ and choose a point $t \in T(K)$. Then there is a selection of points $q_i \in Q_i, i \in K$, and an element $R \in \mathcal{R}$ such that $q_i \in R$ for all $i \in K$, the rows in $\mathcal{D}(R)$ with labels corresponding to $q_i, i \in K$ have a common overlap, and if $r$ is the label of a row to which $q_i$ corresponds then there exits $m_i \in [0, s-1]$ such that $t=r+m_i$.

Let $C$ be the column of $\mathcal{D}(R)$ that is formed by the intersection of $\mathcal{D}(R)$ and the vertical line passing through $t$ on the coordinate axis. Then
\[
C=\{(r, m) \in R \times [0, s-1]: t=r +m\}.
\]

Let $i \in K$. Then $t=r+m_i=q_i+m_i$, hence $(q_i, m_i) \in C$, and so $q_i \in Q_i \cap \theta(C).$ On the other hand, suppose for some $i \in [1, k]$, $q_0 \in Q_i \cap \theta(C).$ Then there is a row of $\mathcal{D}(R)$ labeled by $r_0 \in R$ and $m_0 \in [0, s-1]$ such that $r_0=q_0$ and $(r_0, m_0) \in C$. Hence $t=r_0+m_0=q_0+m_0$, and so $t \in S_i$, whence $i \in K$. We conclude that $K=K(C)$.

Now let $R \in \mathcal{R}$ and let $C$ be a column of $\mathcal{D}(R)$. If $(r, m)$ and $(r^{\prime}, m^{\prime})$ are elements of $C$, then $r+m=r^{\prime}+ m^{\prime}$, hence let $t$ be the common value of $r+m$ as $(r, m)$ varies throughout $C$. If $i \in K(C)$ then there exists $q \in Q_i, r \in \theta(C),$ and $m \in [0, s-1]$ such that $q=r$ and $t=r+m=q+m$, and so $t \in S_i$. Suppose that there exits $i \in [1,  k]$ such that $t \in S_i$. Then there exist $q \in Q_i$ and $m \in [0, s-1]$ such that $t=q+m$. It follows that $q \in R$, and so $(q, m) \in C$. Consequently $q \in Q_i \cap \theta(C)$, i.e., $i \in K(C)$. We conclude that 
\[
t \in \Big(\bigcap_{i\in K(C)}\ S_i \Big)\cap \Big(\bigcap_{i\in [1, k]\setminus K(C)}\ ( \textbf{Q}\setminus S_i)\Big), 
\]
and so $K(C) \in \mathcal{K}_{\max}$.   $\    \Box$ 
\section{The density of $\Pi_+(\textbf{a}, \textbf{b})$}
We begin with the following general situation and then specialize it to the case of interest to us here, namely that of the set $\Pi_+(\textbf{a}, \textbf{b})$. The ability to explain our reasoning concisely in the sequel will be enhanced if we employ the following notation: if $A$ is a nonempty finite set and $ \mathcal{A} \subseteq 2^A$, then $\mathcal{U}(\mathcal{A})$ will denote the set formed by the union of all the elements of $\mathcal{A}$, and $\mathcal{P}(A, 2)$ will denote the set of all 2-block partitions of $A$.

Let $\textbf{S}$ denote a set $\{S_1, \dots, S_m\}$ of nonempty, finite subsets of $[1, +\infty)$ such that for each $i$, all elements of $S_i$ are square-free.. We will say that $p$ is an \emph{allowable prime} (with respect to $\textbf{S})$  if no element of $\bigcup_i S_i$ has $p$ as a factor, and then we will let 
\[
\Pi_+(\textbf{S})
\]
denote the set of all allowable primes such that for all $i$, $S_i$ is either a set of quadratic residues of $p$ or a set of quadratic non-residues of $p$. This is a generalization of the set $\Pi_+(\textbf{a}, \textbf{b})$, as is clear from the following lemma:
\begin{lem}
\label{lem1}$([2, Lemma\   4.1])$ If $(\textbf{a}, \textbf{b})$ is a $2m$-tuple as defined in the begining of section 1 and $\{b_1,\dots,b_k\}$ is the set of distinct values of the coordinates of $\textbf{b}$ then $\Pi_+(\textbf{a}, \textbf{b})$ consists precisely of all primes allowable  with respect to $\{b_1,\dots,b_k\}$ such that  each of the sets $\{b_i: i \in I\}, I \in \Lambda(\mathcal{K})$, is either a set of quadratic residues of $p$ or a set of quadratic non-residues of $p$.
\end{lem}

For each $i \in [1, m]$, let $X_i^+$ (respectively, $X_i^-$) denote the set of all allowable primes $p$  such that $\chi_p$, when restricted to $S_i$, is identically 1, (respectively, is identically $-1$). Let                                                                                                    
\[
M_0=\{i \in [1, m]: 1 \notin S_i\},
\]
\[
M_1=\{i \in [1, m]: 1 \in S_i\}.
\]
Because $X_i^-$ is empty if $i \in M_1$, it follows that
\begin{eqnarray*}
\Pi_+(\textbf{S})&=& \bigcap_{i=1}^m\ \Big(X_i^+ \cup X_i^-\Big)\\
&=&\Big( \bigcap_{i \in M_1}\ X_i^+ \Big) \cap \Big(\bigcap_{i \in M_0}\ \Big(X_i^+ \cup X_i^- \Big) \Big).
\end{eqnarray*}
Let $M_0=\{i_j: j \in [1, \mu]\}$, where $\mu=|M_0|.$ Because $X_i^+$ is disjoint from $X_i^-$, for all $i$, we may write
\[
\bigcap_{i \in M_0}\ \Big(X_i^+ \cup X_i^- \Big)=\bigcup_{(Z_1, \dots, Z_\mu): Z_j \in \{X_{i_j}^+,\ X_{i_j}^-\}, \forall j}\ Z_1 \cap \cdots \cap Z_\mu,
\]
and this union is pairwise disjoint. We rearrange this union so that $\Pi_+(\textbf{S})$ is the pairwise disjoint union of the sets
\begin{equation*}
\bigcap_{1=1}^m\ X_i^+,\tag{3.1}
\]
\begin{equation*}
\Big(\bigcap_{i \in M_1}\ X_i^+ \Big) \cap \Big(\bigcap_{i \in M_0}\ X_i^- \Big), \tag{3.2}
\]
\begin{equation*}
\Big(\bigcap_{i \in M_1\cup P_1}\ X_i^+ \Big) \cap \Big(\bigcap_{i \in P_2}\ X_i^- \Big),\   \tag{3.3}
\]
\begin{equation*}
\Big(\bigcap_{i \in P_1}\ X_i^- \Big) \cap \Big(\bigcap_{i \in M_1 \cup P_2}\ X_i^+ \Big),\   \{P_1, P_2\} \in \mathcal{P}(M_0, 2) \tag{3.4},
\]
In order to calculate the density of  $\Pi_+(\textbf{S})$, it hence suffices to calculate the densities of each of these sets and then add everything up.

We proceed to do precisely that. Observe first that $\bigcap_{i=1}^m X_i^+$  is the set of all allowable primes $p$ such that $\bigcup_{i=1}^m S_i$ is a set of quadratic residues of $p$. Lemma 2.2 hence provides a way to calculate the density of this set. After letting $\pi(z)$ denote the set of prime factors of a square-free integer $z$,  we set
\[
\Pi=\bigcup\  \Big\{\pi(z): z \in \bigcup_i\ S_i \Big\},
\]
\[
n=|\Pi|,\  \textrm{and},
\]
\[
\mathcal{S}_i=\{ \pi(z): z \in S_i \},\  i \in [1, m].
\]
If $F$ is the Galois field of 2 elements and $v: 2^{\Pi} \rightarrow F^n$ is the bijection defined in section 2, then we conclude from Lemma 2.2 that if
\[
d=\textrm{the dimension of the linear span of }\  v\Big( \Big( \bigcup_i\ \mathcal{S}_i \Big)\setminus \{ \emptyset \}\Big)\  \textrm{in}\ F^n
\]
 then
\begin{equation*}
\textrm{the density of}\  \bigcap_{i=1}^m\  X_i^+\  \textrm{is}\  2^{-d}.\tag{3.5}
\]

The next step in our calculation is to compute the density of each set in (3.2)-(3.4). Let $\{P_1, P_2\}$ be a partition of $M_0$ . We first decompose each of these sets into a useful pairwise disjoint union. Toward that end, let
\vspace{0.3cm}
\begin{quote}
$\mathcal{N}(M_0, M_1)=\{N \subseteq \Pi:|N \cap S|\ \textrm{is even}, \textrm{for all}\ S \in \bigcup_{i \in M_1}\ \mathcal{S}_i$
\end{quote}
\begin{center}
$ \textrm{and}\ |N \cap S|\  \textrm{is odd, for all}\  S \in \bigcup_{i \in M_0}\ \mathcal{S}_i\}$,
\end{center}
\vspace{0.1cm}
\begin{quote}
$\mathcal{N}_e(\{P_1, P_2\})= \{N \subseteq \Pi: |N \cap S|\ \textrm{ is even, for all}\  S \in \bigcup_{i \in M_1 \cup P_1} \mathcal{S}_i$
\end{quote}
\begin{center}
$ \textrm{ and} |N \cap S|\  \textrm{is odd, for all}\  S \in \bigcup_{i \in P_2} \mathcal{S}_i\},$ 
\end{center}
\begin{quote}
\vspace{0.1cm}
$\mathcal{N}_o(\{P_1, P_2\})= \{N \subseteq \Pi: |N \cap S|\ \textrm{ is odd, for all}\  S \in \bigcup_{i \in P_1} \mathcal{S}_i$
\end{quote}
\begin{center}
$ \textrm{ and} |N \cap S|\  \textrm{is even, for all}\  S \in \bigcup_{i \in M_1 \cup P_2} \mathcal{S}_i\}.$ 
\end{center}
If for a prime $p$ we set
\[
N(p)= \{q \in \Pi:  \chi_p(q)=-1\},
\]
then
\[
\Big(\bigcap_{i \in M_1}\ X_i^+ \Big) \cap \Big(\bigcap_{i \in M_0}\ X_i^- \Big)=\bigcup_{N \in \mathcal{N}(M_0, M_1)}\  \{p: N(p)=N\},
\]
\[
\Big(\bigcap_{i \in M_1\cup P_1}\ X_i^+ \Big) \cap \Big(\bigcap_{i \in P_2}\ X_i^- \Big)=\bigcup_{N \in \mathcal{N}_e(\{P_1, P_2\})}\  \{p: N(p)=N\},
\]
\[
\Big(\bigcap_{i \in P_1}\ X_i^- \Big) \cap \Big(\bigcap_{i \in M_1 \cup P_2}\ X_i^+ \Big)=\bigcup_{N \in \mathcal{N}_o(\{P_1, P_2 \})}\  \{p: N(p)=N\},
\]
and each of these unions is pairwise disjoint. Now by virtue of Lemma 2.1, the density of each set in these three unions is $2^{-n}$. Observe next that the set in (3.2) is nonempty only if 
\[
\Big( \bigcup_{i\in M_0} S_i \Big) \cap \Big( \bigcup_{i\in M_1}\ S_i \Big)= \emptyset,
\]
and because the elements of the sets $S_i$ are square-free it follows that  this holds if and only if
\[
\Big( \bigcup_{i\in M_0}\ \mathcal{S}_i \Big) \cap \Big( \bigcup_{i\in M_1}\ \mathcal{S}_i \Big)= \emptyset.
\]
Consequently, after observing that the dimension of the linear span of
\[
v \Big( \Big (\Big( \bigcup_{i\in M_0}\ \mathcal{S}_i \Big) \cup \Big(  \bigcup_{i\in M_1}\ \mathcal{S}_i \Big) \Big) \setminus \{\emptyset\} \Big)=v \Big( \Big( \bigcup_i\ \mathcal{S}_i \Big) \setminus \{\emptyset\} \Big)
\]
in $F^n$ is $d$, we conclude from Lemma 2.3 that the cardinality of $\mathcal{N}(M_0, M_1)$  is either $2^{n-d}$ or 0. Following a similar line of reasoning, we find that the cardinality of  $\mathcal{N}_e(\{P_1, P_2\})$ and, respectively, $\mathcal{N}_o(\{P_1, P_2 \}$, is also either $2^{n-d}$ or 0. It follows that
\begin{equation*}
\textrm{the density of}\  \Big(\bigcap_{i \in M_1}\ X_i^+ \Big) \cap \Big(\bigcap_{i \in M_0}\ X_i^- \Big)=\left\{\begin{array}{cc}2^{-d},\ \textrm{if $\mathcal{N}(M_0, M_1) \not= \emptyset$,}\\
0,\ \textrm{otherwise,}\end{array}\right. \tag{3.6}
\]
\begin{equation*}
\textrm{the density of}\  \Big(\bigcap_{i \in M_1 \cup P_1}\ X_i^+ \Big) \cap \Big(\bigcap_{i \in P_2}\ X_i^- \Big)=\left\{\begin{array}{cc}2^{-d},\ \textrm{if $\mathcal{N}_e(\{P_1, P_2\}) \not= \emptyset$,}\\
0,\ \textrm{otherwise,}\end{array}\right. \tag{3.7}
\]
\begin{equation*}
\textrm{the density of}\  \Big(\bigcap_{i \in P_1}\ X_i^- \Big) \cap \Big(\bigcap_{i \in M_1 \cup P_2}\ X_i^+ \Big)=\left\{\begin{array}{cc}2^{-d},\ \textrm{if $\mathcal{N}_o(\{P_1, P_2\}) \not= \emptyset$,}\\
0,\ \textrm{otherwise.}\end{array}\right. \tag{3.8}
\]
Summing the densities from (3.5)-(3.8) now yields the following lemma:
\begin{lem}
\label{lem2}
The density of $\Pi_+(\textbf{S})$ is
\begin{center}
$2^{-d} \Big(1+ \varepsilon+\big| \{ \{P_1, P_2\} \in \mathcal{P}(M_0, 2): \mathcal{N}_e(\{P_1, P_2\}) \not= \emptyset  \} \big|+$
\end{center}
\begin{center}
$\big| \{ \{P_1, P_2\} \in \mathcal{P}(M_0, 2): \mathcal{N}_o(\{P_1, P_2\}) \not= \emptyset  \} \big| \Big),$
\end{center}
where
\[
\varepsilon= \left\{\begin{array}{cc}1,\ \textrm{if $\mathcal{N}(M_0, M_1) \not= \emptyset$,}\\
0,\ \textrm{otherwise.}\end{array}\right.
\]
\end{lem}

At this juncture, the combinatorial parameters that occur in Lemma 3.2 are rather daunting to compute, so as to proceed further, we specialize to the case which is obtained when Lemma 3.2, by way of Lemma 3.1, is applied to a standard $2m$-tuple $(\textbf{a}, \textbf{b})$. The first thing to be done is to produce a set $\textbf{S}$ that can be used in Lemma 3.2 to find the density of $\Pi_+(\textbf{a}, \textbf{b})$. We will proceed as follows: assume to begin with that at least one square-free part of the coordinates of $\textbf{b}$ is not 1; otherwise all of these coordinates are squares, whence the density of $\Pi_+(\textbf{a}, \textbf{b})$ is clearly 1. Let $\sigma_i$ denote the square-free part of $b_i, i \in [1, k]$. For each element $I$ of $\Lambda(\mathcal{K})$, we let $S(I)$ denote the set  formed from the integers $\sigma_i$ for $i \in I$ and then we choose a nonempty subset $Z(I)$ of $I$ such that 
\[
S(I)=\{\sigma_i: i \in Z(I)\}.
\]
 Because $\Pi_+(\textbf{a}, \textbf{b})$ is unaffected by the elements  $I$ of $\Lambda(\mathcal{K})$ for which $S(I)$ is a singleton, if we hence remove all such elements from
$\Lambda(\mathcal{K})$, then for each set $I$ in the set  $\mathcal{I}$ of elements of $\Lambda(\mathcal{K})$  which remain, it follows that $|S(I)| \geq 2$. 

Next, on $\mathcal{I}$ , we define an equivalence relation $\approx$  as follows: if $(I,  J) \in \mathcal{I} \times \mathcal{I}$, then declare that $I \approx J$ if  $S(I)=S(J)$. Select one representative element of  $\mathcal{I}$ from each equivalence class of $\approx$ and let $\Lambda^{\prime}(\mathcal{K})$ denote the set of elements so chosen. The sets $S(I)$ for $I \in \Lambda^{\prime}(\mathcal{K})$ all have cardinality at least 2 and are distinct, and it follows from Lemma 3.1 that $\Pi_+(\textbf{a}, \textbf{b})$ consists precisely of all primes $p$ allowable with respect to $\{b_1, \dots, b_k\}$ such that each of the sets $S(I)$ for $I \in \Lambda^{\prime}(\mathcal{K})$ is either a set of quadratic residues of $p$ or a set of  quadratic non-residues of $p$. Our intension is to use Lemma 3.2 to calculate the density of $ \Pi_+(\textbf{a}, \textbf{b})$ by letting the set $\{ S(I): I \in \Lambda^{\prime}(\mathcal{K}) \}$ play the role of the set $\textbf{S}$.

We now impose the following arithmetic condition on the square-free parts of the coordinates of $\textbf{b}$. Let
\[
\Sigma=\bigcup_{I \in \Lambda^{\prime}(\mathcal{K})}\  Z(I).
\]
We assume that 
\begin{quote}
(3.9) if $\{i, j\} \subseteq \Sigma$ then $\sigma_i \not= \sigma_j$ and the product of all elements in any nonempty subset of $\{\sigma_i: i \in \Sigma\} \setminus \{1\}$ is never a square.
\end{quote}
It is an easy consequence of equation (2.1) in section 2 that this is equivalent to the following combinatorial condition on the prime factors of $\sigma_i, i \in \Sigma$: if $\{i, j\} \subseteq \Sigma$ then $\pi(\sigma_i )\not= \pi(\sigma_j)$ and the repeated symmetric difference of all the elements in every nonempty subset of $\{\pi(\sigma_i): i \in \Sigma\} \setminus \{\emptyset\}$ is never empty.

We proceed to apply Lemma 3.2 with $\textbf{S}=\{ S(I): I \in \Lambda^{\prime}(\mathcal{K}) \}$. In order to do that, we study first the sets which determine the values of the combinatorial parameters in Lemma 3.2. For each $I \in \Lambda^{\prime}(\mathcal{K})$, set 
\[
\mathcal{S}(I)=\{\pi(\sigma): \sigma \in S(I)\}.
\vspace{0.3cm}
\]
Now let
\[
\mathcal{M}_0=\{I \in \Lambda^{\prime}(\mathcal{K}): 1 \notin S(I) \},
\vspace{0.2cm}
\]
\[
\mathcal{M}_1=\{I \in \Lambda^{\prime}(\mathcal{K}): 1 \in  S(I) \}.
\vspace{0.2cm}
\]
N.B. Upon replacing the sets $[1, m]$, $\textbf{S}$,  $M_0$ and $M_1$  defined above by, respectively, $\Lambda^{\prime}(\mathcal{K})$, $\{S(I): I \in \Lambda^{\prime}(\mathcal{K})\}$, $\mathcal{M}_0$, and $\mathcal{M}_1$, we define $\mathcal{N}(\mathcal{M}_0, \mathcal{M}_1),\ \mathcal{N}_e(\{P_1, P_2\})$ and $\mathcal{N}_o(\{P_1, P_2\})$ accordingly. We set
\[
\mathcal{T}(\mathcal{M}_i)=\bigcup_{I \in \mathcal{M}_i}\ \mathcal{S}(I),\ i \in [0, 1],
\]
and suppose that
\begin{equation*}
\mathcal{T}(\mathcal{M}_0) \cap \mathcal{T}(\mathcal{M}_1)= \emptyset. 
\] 
We will use Lemma 2.4 with
\[
A=\bigcup_{i \in \Sigma}\ \pi(\sigma_i),
\vspace{0.2cm}
\]
\[
\mathcal{S}=\mathcal{T}(\mathcal{M}_0),\ \textrm{and}
\vspace{0.2cm}
\]
\[
\mathcal{T}=\mathcal{T}(\mathcal{M}_1),
\]
to prove that
\[
\mathcal{N}(\mathcal{M}_0, \mathcal{M}_1) \not= \emptyset.
\]

Hence, let $U$ be a subset of
\[
\mathcal{T}(\mathcal{M}_0)\cup \mathcal{T}(\mathcal{M}_1) \cup \{\emptyset \}=\big\{\pi(\sigma_i): i \in \Sigma \big\} \cup \{\emptyset\}
 \]
of odd cardinality such that the cardinality of
\[
U \cap \Big( \mathcal{T}(\mathcal{M}_1) \cup \{\emptyset\} \Big)
\]
is even. Then $U \not= \emptyset \not= U \setminus\{\emptyset\} $. We must prove that the repeated symmetric difference of the sets $X \in U$ is nonempty, and  that is an immediate consequence of condition (3.9).

A similar argument also shows that if $\{P_1, P_2\} \in \mathcal{P}(\mathcal{M}_0, 2)$ and
\[
\Big(\bigcup_{I \in \mathcal{M}_1\cup P_1}\ \mathcal{S}(I) \Big) \cap \Big(\bigcup_{I \in P_2}\ \mathcal{S}(I) \Big)=\emptyset
\]
then
\[
\mathcal{N}_e(\{P_1, P_2\}) \not= \emptyset
\]
and if
\[
\Big(\bigcup_{I \in  P_1}\ \mathcal{S}(I) \Big) \cap \Big(\bigcup_{I \in \mathcal{M}_1 \cup P_2}\ \mathcal{S}(I) \Big)=\emptyset
\]
then
\[
\mathcal{N}_o(\{P_1, P_2\}) \not= \emptyset.
\]
Due to the fact that all $\sigma_i$  are square-free, we conclude that $\mathcal{N}(\mathcal{M}_0, \mathcal{M}_1) \not= \emptyset$ if and only if
\[
\Big(\bigcup_{I \in \mathcal{M}_0}\ Z(I) \Big) \cap \Big(\bigcup_{I \in \mathcal{M}_1}\ Z(I) \Big)=\emptyset,
\]
$\mathcal{N}_e(\{P_1, P_2\}) \not= \emptyset$ if and only if
\begin{equation*}
\Big(\bigcup_{I \in \mathcal{M}_1 \cup P_1}\ Z(I) \Big) \cap \Big(\bigcup_{I \in P_2}\ Z(I) \Big)=\emptyset
\tag{3.10},
\]
and $\mathcal{N}_o(\{P_1, P_2\}) \not= \emptyset$ if and only if
\begin{equation*}
\Big(\bigcup_{I \in  P_1}\ Z(I) \Big) \cap \Big(\bigcup_{I \in\mathcal{M}_1 \cup P_2}\ Z(I) \Big)=\emptyset \tag{3.11}.
\]
If we hence let $n=|A|$, $v: 2^A \rightarrow F^n$ be the bijection that is defined in section 2, let $d$ denote the dimension of the linear span in $F^n$ of the set $\{v(\pi(\sigma_i)): i \in \Sigma\} \setminus \{0\}$, and let
\[
\mathcal{P}_{\emptyset}(\mathcal{M}_0, 2)= \big\{ \{P_1, P_2\} \in \mathcal{P}( \mathcal{M}_0, 2): (3.10)\ \textrm{holds for}\ \{P_1, P_2\}\big\},
\]
then Lemma 3.2 implies that
\begin{equation*}
\textrm{the density of}\  \Pi_+(\textbf{a}, \textbf{b})=2^{-d}\big(1+ \varepsilon+2\big|\mathcal{P}_{\emptyset}(\mathcal{M}_0, 2)\big|\big) \tag{3.12},
\]
where
\begin{equation*}
\varepsilon= \left\{\begin{array}{cc}1,\ \textrm{if $\Big(\bigcup_{I \in \mathcal{M}_0}\ Z(I) \Big) \cap \Big(\bigcup_{I \in \mathcal{M}_1}\ Z(I) \Big)= \emptyset$,}\\
0,\ \textrm{otherwise.}\end{array}\right. \tag{3.13}
\]

We turn now to the calculation of the parameters in (3.12). We adopt the following convention, to be used in the rest of what follows; we identify an element $I$ of $\Lambda^{\prime}(\mathcal{K})$ with its subset $Z(I)$ and then, in an abuse of notation that we hope will not prove confusing, continue to denote that element by $I$. As $|S(I)| \geq 2$ for each $I \in \Lambda^{\prime}(\mathcal{K})$, we may assume that $|I| \geq 2$ for all $I \in \Lambda^{\prime}(\mathcal{K})$. We will begin with a calculation of the cardinality of $\mathcal{P}_{\emptyset}(\mathcal{M}_0, 2)$. The set defined by
\[
\mathcal{P}_{\emptyset}(\Lambda^{\prime}(\mathcal{K}), 2)=\Big\{ \{P_1, P_2\} \in \mathcal{P}(\Lambda^{\prime}(\mathcal{K}), 2): \mathcal{U}(P_1) \cap \mathcal{U}(P_2)= \emptyset\Big\}
\]
will be an auxiliary to our analysis and so we will first study the structure of the elements of this set. 

In order to do that, we begin by defining an equivalence relation $\simeq$ on $\Lambda^{\prime}(\mathcal{K})$ as follows: if $X$ and $Y$ are elements of  $\Lambda^{\prime}(\mathcal{K})$ then $X \simeq Y$ if either $X=Y$ or there exists a sequence $\textbf{X}=(X_1, \dots, X_m)$ of elements of $\Lambda^{\prime}(\mathcal{K})$ such that $X_1=X, X_m=Y$ and for all $i \in [1, m-1], X_i \cap X_{i+1} \not= \emptyset$. We will say that the sequence $\textbf{X}$ \emph{connects} $X$ \emph{to} $Y$. Let $\Lambda_1, \dots, \Lambda_{\mu}$ be an enumeration of the equivalence classes of $\simeq$. Then
\begin{equation*}
\Lambda^{\prime}(\mathcal{K})=\bigcup_i\ \Lambda_i, \tag{3.14}
\]
\begin{equation*}
\Sigma= \bigcup_i\ \mathcal{U}(\Lambda_i),\ \textrm{and} \tag{3.15}
\]
\begin{equation*}
\mathcal{U}(\Lambda_i) \cap \mathcal{U}(\Lambda_j)= \emptyset, \textrm{for}\ i \not= j. \tag{3.16}
\]
\begin{lem}
\label{lem3}
The set $\mathcal{P}_{\emptyset}(\Lambda^{\prime}(\mathcal{K}), 2)$ consists precisely of all sets of the form
\[
\Big\{\bigcup_{i \in U_1}\  \Lambda_i,\  \bigcup_{i \in U_2}\  \Lambda_i \Big\}
\]
where $\{U_1, U_2\}$ varies throughout the set $\mathcal{P}([1, \mu], 2)$.
\end{lem}

\emph{Proof}. It follows from (3.14) and (3.16) that every set in the conclusion of Lemma 3.3 is an element $\mathcal{P}_{\emptyset}(\Lambda^{\prime}(\mathcal{K}), 2)$.

Let $\{P_1, P_2\} \in \mathcal{P}_{\emptyset}(\Lambda^{\prime}(\mathcal{K}), 2)$. We claim that for each $n \in [1, \mu]$, either $P_1 \cap \Lambda_n$ or $P_2 \cap \Lambda_n$ is empty.
Suppose not and find an $n \in[1, \mu]$ such that
\[
L_i=P_i \cap \Lambda_n \not= \emptyset,\ \textrm{for}\ i \in [1, 2].
\]
Then $\{L_1, L_2\}$ is a partition of $\Lambda_n$ such that
\begin{equation*}
\mathcal{U}(L_1) \cap \mathcal{U}(L_2)= \emptyset. \tag{3.17}
\]

Let $X \in L_1, Y \in \Lambda_n$. We wish to prove that $Y \in L_1$, and so we may assume that there exists a sequence $(X_1, \dots, X_m)$ connecting $X$ to $Y$. Suppose inductively that $j>1$ and $X_j \in L_1$. Let $i \in X_j \cap X_{j+1}$. Then $i \in \mathcal{U}(L_1)$, and if $X_{j+1} \in L_2$ then $i \in \mathcal{U}(L_2)$, which contradicts (3.17). Hence $X_{j+1} \in L_1$, and we conclude by induction that $X_j \in L_1$, for all $j \in [1, m]$. In particular, $Y=X_m \in L_1$. However, $Y$ is an arbitrary element of $\Lambda_n$, and so it follows that $\mathcal{U}(L_1)=\mathcal{U}(\Lambda_n)$, contrary to the hypothesis that $L_2$ is not empty. This verifies our claim.

 It follows that if $i \in [1, 2], n \in [1, \mu],$ and $P_i \cap \Lambda_n\not= \emptyset$, then $\Lambda_n \subseteq P_i$; otherwise there exists $n \in [1, \mu]$ such that $P_1 \cap  \Lambda_n \not= \emptyset \not=P_2 \cap  \Lambda_n$, contradicting what we have just shown. Because of (3.14), we hence conclude that for $i\in [1, 2]$, $\{n \in [1, \mu]: P_i \cap  \Lambda_n \not= \emptyset \} \not= \emptyset$  and $P_i$ is the union of the $\Lambda_n$'s with which it has a nonempty intersection. This verifies the lemma.  $\    \Box$
\vspace{0.3cm}

We can now calculate the cardinality of $\mathcal{P}_{\emptyset}(\mathcal{M}_0, 2)$. Suppose first that $\mathcal{M}_1=\emptyset$. Then  $\mathcal{P}_{\emptyset}(\mathcal{M}_0, 2)=\mathcal{P}_{\emptyset}(\Lambda^{\prime}(\mathcal{K}), 2)$, and so we deduce from Lemma 3.3 that
\begin{equation*}
\big|\mathcal{P}_{\emptyset}(\mathcal{M}_0, 2)\big|=2^{\mu-1}-1,\ \textrm{if}\  \mathcal{M}_1=\emptyset. \tag{3.18}
\]

Suppose next that $\mathcal{M}_1\not=\emptyset,$ i.e., there is a unique element $i_0$ of $\Sigma$ such that $\sigma_{i_0}=1$. Because of (3.15) and (3.16) there is a unique $n_0 \in [1, \mu]$ such that $i_0 \in \Lambda_{n_0}$, and so
\begin{equation*}
\mathcal{M}_1= \{I \in \Lambda_{n_0}: i_0 \in I \}, \tag{3.19}
\]
\begin{equation*}
\mathcal{M}_0=  \Big( \Lambda_{n_0} \setminus \mathcal{M}_1 \Big) \cup \Big( \bigcup_{n\not= n_0}\  \Lambda_n \Big). \tag{3.20}
\]

Assume first that $\Lambda_{n_0} \setminus \mathcal{M}_1 \not= \emptyset$. Let $\{P_1, P_2\} \in \mathcal{P}_{\emptyset}( \mathcal{M}_0, 2)\}$. Then
\[
\{\mathcal{M}_1 \cup P_1, P_2\} \in \mathcal{P}_{\emptyset}(\Lambda^{\prime}(\mathcal{K}), 2),
\]
and so by Lemma 3.3, there exists a partition $\{U_1, U_2\}$ of $[1, \mu]$ such that
\[
\mathcal{M}_1 \cup P_1= \bigcup_{i \in U_1}\  \Lambda_i,\  \textrm{and}
\]
\[
P_2=\bigcup_{i \in U_2}\  \Lambda_i.
\]
It follows that $n_0 \in U_1$ and
\[
P_1=\Big(\Lambda_{n_0} \setminus \mathcal{M}_1\Big) \cup \Big( \bigcup_{i \in U_1 \setminus \{n_0\}}\   \Lambda_i\Big).
\]
Because $\Lambda_{n_0} \setminus \mathcal{M}_1 \not=\emptyset$, it follows that there is a bijection of $\mathcal{P}_{\emptyset}( \mathcal{M}_0, 2)$ onto  $\mathcal{P}([1, \mu], 2)$. Hence
\begin{equation*}
\big|\mathcal{P}_{\emptyset}( \mathcal{M}_0, 2)\big|=2^{\mu-1}-1,\ \textrm{if}\ \emptyset \not= \mathcal{M}_1 \not= \Lambda_{n_0}.\tag{3.21}
\]

Assume now that $\mathcal{M}_1=\Lambda_{n_0}$. Applying Lemma 3.3 to a partition $\{P_1, P_2\} \in \mathcal{P}_{\emptyset}(\mathcal{M}_0, 2)$ as before, we find a proper subset $U$ of $[1, \mu]$ such that
\[
n_0 \in U\  \textrm{and}\  P_1=\bigcup_{i \in U \setminus \{n_0\}}\ \Lambda_i.
\]
Because $P_1$ is nonempty, it follows that $ U \setminus \{n_0\}$ is also nonempty. We conclude that there is a bijection of $\mathcal{P}_{\emptyset}(\mathcal{M}_0, 2)$ onto $\mathcal{P}([1, \mu-1], 2)$, and so
\begin{equation*}
\big|\mathcal{P}_{\emptyset}(\mathcal{M}_0, 2)\big|=2^{\mu-2}-1,\ \textrm{if}\ \emptyset \not=\mathcal{M}_1= \Lambda_{n_0}. \tag{3.22}
\]

We turn next to the calculation of the parameter $\varepsilon$ in (3.12). As a consequence of (3.13),
\begin{equation*}
\varepsilon=1\ \textrm{if}\ \mathcal{M}_1= \emptyset. \tag{3.23}
\]

Suppose that $\mathcal{M}_1\not= \emptyset.$ In order to compute $\varepsilon$ in this case, we must determine when 
\[
\mathcal{U}(\mathcal{M}_0) \cap \mathcal{U}(\mathcal{M}_1)= \emptyset.
\]
This is accomplished in the following lemma.
\begin{lem}
\label{lem4}
$\mathcal{U}(\mathcal{M}_0)$ is disjoint from $ \mathcal{U}(\mathcal{M}_1)$ if and only if
\begin{equation*}
\mathcal{M}_1=\Lambda_{n_0}\ \textrm{and}\  \mathcal{M}_0= \bigcup_{n\not= n_0}\  \Lambda_n. \tag{3.24}
\]
\end{lem}

\emph{Proof}. That $\mathcal{U}(\mathcal{M}_0)$ is disjoint from $\mathcal{U}(\mathcal{M}_1)$ when (3.24) is true follows immediately from (3.16). On the other hand, if  $\mathcal{U}(\mathcal{M}_0)$ and $\mathcal{U}(\mathcal{M}_1)$ are disjoint, then 
\[
\mathcal{U}(\mathcal{M}_1) \cap \mathcal{U}(\Lambda_{n_0} \setminus \mathcal{M}_1)= \emptyset,
\]
and so if $\Lambda_{n_0} \setminus \mathcal{M}_1\not= \emptyset$, then $\{\mathcal{M}_1, \Lambda_{n_0} \setminus \mathcal{M}_1 \}$ is a partition of $\Lambda_{n_0}$ whose existence was shown to be impossible in the proof of Lemma 3.3. Hence
\[
 \mathcal{M}_1=\Lambda_{n_0} ,
\]
and (3.24) is now a consequence of this equation and (3.20). $\  \Box$
\vspace{0.3cm}
 
 We conclude from (3.23) and Lemma 3.4 that
 \begin{equation*}
 \varepsilon= \left\{\begin{array}{cc}1,\ \textrm{if either $\mathcal{M}_1= \emptyset$ or $\emptyset \not= \mathcal{M}_1= \Lambda_{n_0}$,}\\
0,\ \textrm{if}\ \emptyset \not= \mathcal{M}_1 \not= \Lambda_{n_0}.\end{array}\right. \tag{3.25}
\]

Finally, we evaluate the parameter $d$. Letting $F= \textrm{the Galois field of cardinality}\ 2, A=\bigcup_{i \in \Sigma} \pi(\sigma_i)$ and $n=|A|$ as before,  we first recall that the \emph{support} supp$(v)$ of a vector $v=(v(1),\dots, v(n)) \in F^n$ is the set $\{i \in [1, n]: v(i)=1\}$ and then observe that a nonempty subset $V$ of $F^n \setminus\{0\}$ is linearly independent over $F$ if and only if for each nonempty subset $S$ of $V$, the repeated symmetric difference of the sets supp$(v), v \in S,$ is not empty. If the elements $a_1< \cdots < a_n$ of $A$ are arranged in increasing order and if $\eta: A\rightarrow [1, n]$ is the bijection defined by $\eta(a_i)=i, i \in [1, n]$, then
\[
\textrm{supp}(v(\pi(\sigma_i)))= \eta(\pi(\sigma_i)),\  i \in \Sigma,
\]
and so, as a consequence of assumption (3.9), for each nonempty subset $S$ of the set $\{\eta(\pi(\sigma_i)): i \in \Sigma \} \setminus \{\emptyset\}$, the repeated symmetric difference of the elements of $S$ is not empty. It follows that $\{v(\pi(\sigma_i)): i \in \Sigma\} \setminus \{0\}$ is linearly independent over $F$. When we now observe that $\{v(\pi(\sigma_i)): i \in \Sigma\}$ contains 0 if and only if $\mathcal{M}_1 \not= \emptyset$, if follows that
\begin{equation*}
d= \left\{\begin{array}{cc}|\Sigma|,\ \textrm{if $\mathcal{M}_1= \emptyset$ ,}\\
|\Sigma|-1,\ \textrm{if $\mathcal{M}_1\not= \emptyset$.}\end{array}\right. \tag{3.26}
\]

All of the parameters that are required for the calculation of the density of $\Pi_+(\textbf{a}, \textbf{b})$ by means of equation (3.12) have now been determined, and we see in particular that $d$ and $|\mathcal{P}_{\emptyset}(\mathcal{M}_0, 2)|$ are certain functions of $|\Sigma|$ and the number $\mu$ of equivalence classes of the equivalence relation $\simeq$ defined on $\Lambda^{\prime}(\mathcal{K})$. Although $|\Sigma|$ is obviously determined directly by $\Sigma$, $\mu$ is instead determined directly by an equivalence relation defined on a certain set of subsets of $\Sigma$. It would be preferable if $\mu$ could be calculated by a procedure that ties it more closely to $\Sigma$ itself, and a good way to do that would be to determine $\mu$ by means of an equivalence relation defined directly on $\Sigma$, instead of an equivalence relation defined on a set of subsets of $\Sigma$. Toward that end, we will now define  and study an equivalence relation $\sim$ on $\Sigma$ such that $\mu$ is also equal to the number of equivalence classes of $\sim$. 

In order to define this equivalence relation, it will be necessary to consider a certain subset of columns in the set $\mathcal{C}$ defined by equation (2.2) in section 2. If $C \in \mathcal{C}$ then we let $S(C)$ denote the set formed by the integers $\sigma_i, i \in K(C)$ and declare that $C$ is an \emph{essential column} if $|S(C)| \geq 2$. If $|S(C)|=1$ then we say that $C$ is a \emph{non-essential column}. If $i$ and $j$ are elements of $\Sigma$ then we set $i \sim j$ if either $i=j$ or there exists a sequence $\textbf{s}=(s_1,\dots,s_m)$ of elements of $\Sigma$ such that $s_1=i, s_m=j$ and for all $l \in [1, m-1], s_l \not= s_{l+1}$ and there exits an essential column $C_l$  such  that $\{s_l, s_{l+1}\} \subseteq K(C_l)$. We will say that the sequence $\textbf{s}$ \emph{connects} $i$ \emph{to} $j$ and we will call the coordinates of the sequence $(K(C_1),\dots,K(C_{m-1}))$ the \emph{links of} $\textbf{s}$. We will prove the following proposition.
\begin{prp}
\label{prp5}
Let $\Sigma/ \sim$ and $\Lambda^{\prime}(\mathcal{K})/ \simeq$ denote, respectively, the set of equivalence classes of $\sim$ and $\simeq$. For $\varpi \in \Sigma/ \sim$ let $\Phi(\varpi)$ denote the set
\[
\big\{X \in \Lambda^{\prime}(\mathcal{K}): X \subseteq \varpi \big\}.
\]
Then $\Phi$ defines a bijection of $\Sigma/ \sim$ onto $\Lambda^{\prime}(\mathcal{K})/ \simeq$ whose inverse mapping sends an element $F$ of $\Lambda^{\prime}(\mathcal{K})/ \simeq$ to the union of the elements of $F$.
\end{prp}

It is an immediate consequence of Proposition 3.5 that
\begin{equation*}
\mu=|\Sigma/ \sim|. \tag{3.27}
\]
Moreover, if $i_0$ is the element of $\Sigma$ that determines $\mathcal{M}_1$ as per (3.19), $\varpi_0$ is the equivalence class of $\sim$ which contains $i_0$, and $n_0$ is the element of $[1, \mu]$ such that $i_0 \in \Lambda_{n_0}$, then
\begin{equation*}
\Phi(\varpi_0)= \Lambda_{n_0}. \tag{3.28}
\]

The proof of Proposition 3.5 will be facilitated by use of the following lemma. For the purposes of our arguments it will be convenient to observe that $i \sim j$ if and only if $i=j$ or there exists a sequence $(s_1,\dots,s_m)$ of elements of $\Sigma$ such that $s_1=i, s_m=j$ and for all $l \in[1, m-1]$, there is an element $X_l$ of $\Lambda^{\prime}(\mathcal{K})$ such that $s_l \in X_l$ and $s_{l+1} \in X_l$. We retain the  terminology in this situation, appropriately modified, that was introduced in the definition of $\sim$. 

\begin{lem}
\label{lem6}
$(i)$ If $i \sim j$ and $\textbf{s}$ is a sequence connecting $i$ to $j$ then all the coordinates of any sequence of links of $\textbf{s}$ are contained in the same equivalence class of $\simeq$.

$(ii)$ If $X \simeq Y$ then the union of all of the coordinates of any sequence connecting $X$ to $Y$ is contained in an equivalence class of $\sim$.
\end{lem}
 
 \emph{Proof}. $(i)$ Let $(X_1,\dots, X_{m-1})$ be a sequence of links of $\textbf{s}=(s_1,\dots,s_m)$. If $m=2$ then there is only one link , it is contained in $\Lambda^{\prime}(\mathcal{K})$, hence it is clearly contained in an equivalence class of $\simeq$. If $m>2$ then $s_{i+1} \in X_i \cap X_{i+1}$ for all $i \in [1, m-2]$, and so $X_1,\dots, X_{m-1}$ are all contained in the same equivalence class of $\simeq$.
 
 $(ii)$. Let $(X_1,\dots, X_m)$ be a sequence connecting $X$ to $Y$, let $\{s, s^{\prime}\} \subseteq \bigcup_i X_i$, and choose $j$ and $j^{\prime}$ such that $s \in X_j$ and $s^{\prime} \in X_{j^{\prime}}$. We must prove that $s \sim s^{\prime}$, and if $j=j^{\prime}$ then $(s, s^{\prime})$ is a  sequence connecting $s$ to $s^{\prime}$, hence this is true. Otherwise, assume that $j<j^{\prime}$, say, and choose $s_i \in X_i \cap X_{i+1}$ for each $i \in [j, j^{\prime}-1]$. Then $(X_j,\dots,X_{j^{\prime}})$ is a sequence of links for the sequence $(s, s_j,\dots,s_{j^{\prime}-1}, s^{\prime})$ connecting $s$ to $s^{\prime}$, and so $s \sim s^{\prime}$. $\     \Box$ 
\vspace{0.3cm}

\emph{Proof of Proposition} 3.5. Let $\varpi \in \Sigma/\sim$. In order to show that $\Phi(\varpi) \in \Lambda^{\prime}(\mathcal{K})/ \simeq$, we must prove that $\Phi(\varpi)$ is nonempty, contained in an equivalence class of $\simeq$, and is maximal with respect to this latter property. Toward that end, we take $i \in \varpi$ and then find $X \in \Lambda^{\prime}(\mathcal{K})$ such that $i \in X$. If $\{i, j\} \subseteq X$, then $(i, j)$ connects $i$ to $j$, hence $j \in \varpi$ and so $X \subseteq \varpi$, hence $X \in \Phi(\varpi)$. Suppose next that $\{X, Y\} \subseteq \Phi(\varpi)$, and then choose $i \in X, j \in Y$. Because $i \sim j$, either $i=j$, in which case $i \in X \cap Y$, and so $X \simeq Y$, or there is  a sequence $(s_1,\dots,s_m)$ connecting $i$ to $j$. Because $i$ and $j$ are, respectively, contained in the intersection of $X$ with the first link of this sequence and the intersection of $Y$ with the last link, it follows from Lemma 3.6 $(i)$ that $X \simeq Y$. Finally let $X \in \Phi(\varpi)$, assume that $X \simeq Z \in \Lambda^{\prime}(\mathcal{K})$, and hence let $(X_1,\dots,X_m)$ connect $X$ to $Z$. By virtue of Lemma 3.6$(ii)$, $\bigcup_iX_i$ is contained in an equivalence class $\varpi^{\prime}$ of $\sim$, hence $X \subseteq \varpi \cap \varpi^{\prime}$ and so $ \varpi=\varpi^{\prime}$. But then $Z=X_m \subseteq \varpi^{\prime}= \varpi$, hence $Z \in \Phi(\varpi)$. 

If $\varpi$ and $\varpi^{\prime}$ are elements of $\Sigma/\sim$ and $X \in \Phi(\varpi) \cap \Phi(\varpi^{\prime})$ then $X \subseteq \varpi \cap \varpi^{\prime}$, hence $\varpi=\varpi^{\prime}$. It follows that $\Phi$ is injective.

Let $F \in \Lambda^{\prime}(\mathcal{K})/\simeq$ and set $\varpi=$ the union of the elements of $F$. That $F \subseteq \Phi(\varpi)$ is clear from the definition of $\Phi(\varpi)$. In order to verify the reverse inclusion, let $X \in \Phi(\varpi)$ and choose $i \in X$. There exits a $Y \in F$ such that $i \in Y$, hence $i \in X \cap Y$, and so $X \simeq Y$, hence $X \in F$. It follows that $\Phi(\varpi)=F$, and so to conclude that $\Phi$ is surjective, we need only prove that $\varpi \in \Sigma/\sim$, i.e., we must show that $\varpi$ is not empty, which is clear from its definition, that it is contained in an equivalence class of $\sim$, and that it is maximal with respect to this latter property. We hence take $\{i, j\} \subseteq \varpi$ and choose elements $X$ and $Y$ of $F$ such that $i \in X$ and $j \in Y$. Then $X \simeq Y$ and so, in light of Lemma 3.6$(ii)$, $X \cup Y$ is contained in an equivalence class of $\sim$, hence $i \sim j$. Suppose next that $i \in \varpi$ and $i \sim z \in \Sigma$. Choose $X \in F$ with $i \in X$. Then $X$ intersects the first link in a sequence connecting $i$ to $z$, and so it follows from Lemma 3.6$(i)$ that $X$ is $\simeq$-equivalent to the link $Y$ containing $z$. Thus $Y \in F$, hence $z \in \varpi$. $\    \Box$  
\vspace{0.3cm}

We now deduce from (3.12), (3.18), (3.19), (3.21), (3.22) and (3.25)-(3.28) the following theorem, the principal result of this paper.

\begin{thm}
\label{thm7}
Let $(\mathbf{a}, \mathbf{b})$ be a standard $2m$-tuple and assume that the square-free parts of the coordinates of $\mathbf{b}$ satisfy $(3.9)$. Let
\[
\mathcal{M}_1=\{I \in \Lambda^{\prime}(\mathcal{K}): 1 \in S(I) \},
\vspace{0.2cm}
\]
\[
\sigma=|\Sigma|,\   \textrm{and}\ \mu=\big|\Sigma/ \sim \big|.
\vspace{0.2cm}
\]
If whenever $\mathcal{M}_1$ is nonempty,  $\varpi_0$ is the equivalence class of $\Sigma/ \sim$ that contains the index $i_0$ which determines $\mathcal{M}_1$ as per $(3.19)$, then the density of $\Pi_+(\mathbf{a}, \mathbf{b})$ is
\[
2^{\mu- \sigma},\  \textrm{if $\mathcal{M}_1= \emptyset$ or  $\mathcal{M}_1= \Phi(\varpi_0)$,}
\vspace{0.2cm}
\]
or
\[
2^{1- \sigma}(2^{\mu}-1),\   \textrm{if  $\emptyset \not= \mathcal{M}_1\not= \Phi(\varpi_0)$.}
\]
\end{thm}

Theorem 3.7 shows that  each element of $\Sigma$ contributes a factor of $1/2$ to the density of $\Pi_+(\mathbf{a}, \mathbf{b})$ and each equivalence class of $\sim$ contributes essentially a factor of 2 to the density. If $\mu=0$ then $\Sigma= \emptyset$ , i.e., for all $I \in \Lambda(\mathcal{K})$, $|S(I)|=1$. It follows that $\Pi_+(\mathbf{a}, \mathbf{b})$ is the set of all allowable primes, with density 1. We note incidentally that $|S(I)|=1$ for all $I \in \Lambda(\mathcal{K})$ if and only if $\prod_{i \in I} b_i$ is a square for all $I \in \Lambda(\mathcal{K})$, and this is valid, as we reported in the introduction, if and only if the cardinality of the set $\{A\in AP(\textbf{a}, \textbf{b}; s)\cap 2^{[1, p-1]}: \chi_p(a)=\varepsilon, \ \textrm{for all}\ a\in A\}$ is asymptotic to $(b\cdot 2^{\kappa})^{-1}p$ as $p \rightarrow +\infty$, where $ \varepsilon \in \{-1, 1\}$, $b$ is the largest value of the coordinates of $\textbf{b}$,  and $\kappa$ is the cardinality of $\bigcup_{i=1}^k  S_i $. Because the cardinality of each equivalence class of $\sim$ is at least 2, it follows that whenever $\mu \geq 1$ then $|\Sigma| \geq 2\mu$ and so the density of  $\Pi_+(\mathbf{a}, \mathbf{b})$is at most $2^{-\mu}$ if $\mathcal{M}_1= \emptyset$ or $\emptyset \not= \mathcal{M}_1= \Phi(\varpi_0)$, and at most $(2^{\mu}-1)/2^{2\mu-1}$, otherwise.

The formulae in Theorem 3.7 reduce the calculation of the density of $\Pi_+(\textbf{a}, \textbf{b})$ to the calculation of the the equivalence classes of $\sim$. In general this calculation can be rather complicated; it can be shown that if $k \in [2, + \infty), (Q_1,\dots,Q_k)$ is any fixed  $k$-tuple of nonempty, finite subsets of the non-negative rationals, and $m=\sum_{i=1}^k |Q_i|$ then there exits infinitely many standard $2m$-tuples $(\textbf{a}, \textbf{b})$ such that $\textbf{b}$ has $k$ distinct coordinates $b_1,\dots, b_k$, the square-free parts $\sigma_i$ of $b_i, i \in [1, k]$ have the property that $\pi(\sigma_i)$ is a proper subset of $\pi(\sigma_{i+1})$ for $i \in [1, k-1]$, and $\{a/b_i: a \in A(b_i)\}=Q_i$, for all $i \in [1, k]$. However there is an interesting class of standard $2m$-tuples for which is available a geometric procedure that determines these equivalence classes in a fairly elegant and straightforward manner. We will discuss that next.

Recall from the introduction that if  $k \in [2, + \infty)$ then a standard $2k$-tuple $(\textbf{a}, \textbf{b})$ is said to be admissible if the coordinates of $\textbf{b}$ are distinct and $a_ib_j-a_jb_i \not= 0$ for $i \not= j$. Let $(\textbf{a}, \textbf{b})$ denote a fixed admissible $2k$-tuple. We will assume with no loss of generality that the coordinates of $\textbf{a}$ and $\textbf{b}$ are indexed so that the rational numbers $q_i=a_i/b_i$ are strictly increasing as the index $i$ increases from $1$ to $k$. Because $(\textbf{a}, \textbf{b})$ is admissible, the sets $S_i$ given by (1.1) in the introduction simplify to  
\begin{equation*}
S_i=\{q_i+j: j \in [0, s-1]\}, i \in [1, k],
\]
i.e., in the notation that was used in section 2, we have that $A(b_i)=Q_i= \{q_i\}, i \in [1, k]$ and $Q=\{q_1,\dots, q_k\}$. This leads to the following simplified determination of the the quotient diagram of $(\textbf{a}, \textbf{b})$ and the set $\Lambda(\mathcal{K})$.

Let $\mathcal{R}$ denote the partition of $Q$ that was constructed in section 2. If $R \in \mathcal{R}$ and $q_j$ labels row $i$ of the the overlap diagram $\mathcal{D}(R)$ then we replace the label $(q_j, l)$ of the $l$-th point in row $i$  by $(j, l)$, and then observe that if $C$ is a column of $\mathcal{D}(R)$ and $\theta$ now denotes the canonical projection of $[1, k] \times [0, s-1]$ onto $[1, k]$, then $K(C)=\theta(C)$. As will become clear, it is precisely this fact that will allow us to directly apply the geometry of overlap diagrams to the determination of the equivalence classes of $\sim$. When $(\textbf{a}, \textbf{b})$  is no longer admissible, the set $S(q)=\{i \in [1, k]: q \in Q_i\}$ that each element $q$ of $\theta(C)$ contributes to  $K(C)$ may have cardinality greater than 1, $S(q) \cap S(q^{\prime})$ may be nonempty for particular subsets $\{q, q^{\prime}\}$ of $\theta(C)$, and the map $q \rightarrow S(q)$ may no longer be injective, which facts can significantly complicate the calculation of the equivalence classes. It follows that if $\{R_1,\dots, R_m\}$ is an enumeration of the elements $R$ of $\mathcal{R}$ such that $|R| \geq 2$ then the quotient diagram of $(\textbf{a}, \textbf{b})$ is given by the $m$-tuple of overlap diagrams $(\mathcal{D}(R_1),\dots, \mathcal{D}(R_m))$. We also deduce from Lemma 2.5 that if $\mathcal{C}_n$ denotes the set of all columns of $\mathcal{D}(R_n), n \in [1, m]$ and if we let $\mathcal{C}= \bigcup_n \mathcal{C}_n$ then
\begin{equation*}
\Lambda(\mathcal{K})=\bigcup_{C \in \mathcal{C}}\ \mathcal{E}(\theta(C)). \tag{3.29}
\]

Let $v\in [1, +\infty)$ and for each $n\in [1, v]$, let $\mathcal{D}(n)$ be a fixed but arbitrary overlap diagram with $k_n$ rows, $k_n\geq 2$, and gap sequence $(d(i, n): i\in [1, k_n-1])$, with no gap exceeding $s-1$. Let $k_0=0, k= \sum_n k_n$. We will now exhibit infinitely many admissible $2k$-tuples $(\textbf{a}, \textbf{b})$ whose quotient diagram is $\Delta=(\mathcal{D}(n): n\in [1, v])$. This is done by taking the $(k-1)$-tuple $(d_1,\dots, d_{k-1})$ in the following lemma to be 
\[
d_i=\left\{\begin{array}{cc}d\Big(i-\sum_0^n k_j, n+1\Big),\ \textrm{if}\ n\in [0, v-1]\ \textrm{and}\  i\in \Big[1+\sum_0^n k_j, -1+\sum_0^{n+1} k_j\Big],\\
s,\ \textrm{elsewhere},\end{array}\right.
\]
and then letting $(\textbf{a}, \textbf{b})$ be any $2k$-tuple obtained from the construction in the lemma.

\begin{lem}
\label{lem8}
For $k\in [2, +\infty)$, let $(d_1,\dots, d_{k-1})$ be a $(k-1)$-tuple of positive integers. Define $k$-tuples $(a_1,\dots, a_k), (b_1,\dots, b_k)$ of positive integers inductively as follows: let $(a_1, b_1)$ be arbitrary, and if $i>1$ and $(a_i, b_i)$ has been defined, choose $t_i\in [2, +\infty)$ and set
\[
a_{i+1}=t_i(a_i+d_ib_i),\ \ b_{i+1}=t_ib_i.
\]
Then
\[
a_ib_j-a_jb_i=\Big(\sum_{r=j}^{i-1}\ d_r\Big)b_ib_j,\ \textrm{for all}\  i>j.
\]
\end{lem}

If one chooses $b_1$ and all subsequent $t_i$'s in the construction of Lemma 3.8 to be distinct primes, then one obtains infinitely many admissible $2k$-tuples $(\textbf{a}, \textbf{b})$ with given quotient diagram $\Delta$ such that $\Pi_+(\textbf{a}, \textbf{b})$ and $\Pi_-(\textbf{a}, \textbf{b})$ are both infinite. Moreover, in that construction $b_i$  is square-free for all $i \in [1,k]$ and $\pi(b_i)$ is a proper subset of $\pi(b_{i+1})$, for all $i \in [1, k-1]$, from whence it follows that condition (3.9) is automatically satisfied.This shows that there are infinitely many admissible $2k$-tuples with a given quotient diagram which satisfy all of the hypotheses of Theorem 3.7. In particular, the number of blocks in the quotient diagram of these $2k$-tuples can be arbitrary, and, as we shall see, this number, in certain situations, is an important parameter in the calculation of the density of $\Pi_+(\textbf{a}, \textbf{b})$.

Let $(\textbf{a}, \textbf{b})$ be an admissible $2k$-tuple with quotient diagram $(\mathcal{D}(R_1),\dots, \mathcal{D}(R_m))$ and let $\Lambda^{\prime}(\mathcal{K})$ and $\Sigma$ be the sets as defined above by $(\textbf{a}, \textbf{b})$. In order to study the associated equivalence relation $\sim$ on $\Sigma$, we begin with the following lemma. Its statement will be more concise if we let $D_n$ denote the subset of $[1, k]$ such that $R_n=\{q_i: i \in D_n\}$ and then let $\Sigma_n=\Sigma \cap D_n, n \in [1, m]$.  

\begin{lem}
\label{lem9}
If $E$ is an equivalence class of $\sim$ then there exits $n \in [1, m]$ such that $E \subseteq \Sigma_n$.
\end{lem}

\emph{Proof}. Note first that 
\begin{equation*}
D_n=\bigcup_{C \in \mathcal{C}_n} \theta(C), n\in [1, m], \tag{3.30}
\]
hence it follows from (3.29) that
\[
\mathcal{U}(\Lambda(\mathcal{K}))=\bigcup_n D_n,
\]
from which we conclude that
\begin{equation*}
\Sigma=\bigcup_n \Sigma_n. \tag{3.31}
\]

Now let $i \in E$ and conclude from (3.31) that there is an $n \in [1, m]$ for which $i \in \Sigma_n$. Let $j \in E \setminus \{i\}$ and choose a sequence $(s_1,\dots,s_t)$ connecting $i$ to $j$, with associated link sequence $(\theta(C_1),\dots, \theta(C_{t-1}))$. Because $i \in D_n \cap \theta(C_1)$, from (3.30) and the fact that $D_1,\dots, D_m$ are pairwise disjoint it follows that $\theta(C_1) \subseteq D_n$, hence $s_2 \in D_n$. An induction argument using the fact that $\{s_l, s_{l+1}\} \subseteq \theta(C_l)$ for $l \in [1, t-1]$ now shows that $j=s_t \in \theta(C_{t-1}) \subseteq D_n$. Hence $E \subseteq D_n \cap \Sigma= \Sigma_n$. $\     \Box$
\vspace{0.3cm}

By virtue of Lemma 3.9, in order to determine the equivalence classes of $\sim$, we need only study the action of $\sim$ on each of the subsets $\Sigma_n$ of $\Sigma, n \in [1, m]$. It is here that some simple geometry of the blocks in the quotient diagram of $(\textbf{a}, \textbf{b})$ will be brought to bear. In order to do that, we must first study the essential columns in the blocks of the quotient diagram.
\begin{lem}
\label{lem10}
If $C$ is an essential column in $\mathcal{C}_n$ then $\theta(C)$ is contained in an equivalence class of $\sim$.
\end{lem}

\emph{Proof}. Let $Z(C)$ denote a subset of $\theta(C)$ such that $S(C)=\{\sigma_i: i \in Z(C)\}$. Because the square-free parts $\sigma_i$ for $i \in Z(C)$ are distinct, $|Z(C)| \geq 2$, and $Z(C) \subseteq \theta(C)$, it follows immediately from the definition of $\sim$ that $Z(C)$ is contained in an equivalence class $E$ of $\sim$. If $i \in \theta(C) \setminus Z(C)$ then there exits $j \in Z(C)$ such that $\sigma_i= \sigma_j$, and so if we choose $l \in Z(C)$ such that $l \not= j$ then $\sigma_i \not= \sigma_l$, hence $\{i, l\} \subseteq \Sigma \cap \theta(C)$, from whence it follows that $i \sim l$. We conclude that $\theta(C) \subseteq E$. $\     \Box$
 \vspace{0.3cm}

A particular consequence of Lemma 3.10 is that whenever $C \in \mathcal{C}_n$ and $|S(C)|\geq 2$, then $\theta(C) \subseteq \Sigma_n$. On the other hand, if $X \in \Lambda^{\prime}(\mathcal{K})$ and $X \subseteq D_n$ then there exits $C \in \mathcal{C}_n$ and $I \in \mathcal{E}(\theta(C))$ such that $|S(I)| \geq 2$ and $X=Z(I)$. Hence $|S(C)| \geq 2$, and so it follows that
\begin{equation*}
\Sigma_n=\bigcup_{C \in \mathcal{C}_n: |S(C)| \geq 2} \theta(C). \tag{3.32}
\]
It is worthwhile to note that it also follows from the reasoning which produced (3.32) that if, for a standard $2m$-tuple, $\mathcal{C}$ is the set of columns defined by (2.2) in section 2,  then, in general,
\[
\Sigma=\bigcup_{C \in \mathcal{C}: |S(C)| \geq 2} K(C).
\] 

From $\mathcal{D}(R_n)$, we now remove all of the non-essential columns, i.e. the columns $C$ such that $|S(C)|=1$. If the subset of points which remain in $\mathcal{D}(R_n)$ is not empty then we let $\mathcal{D}_n ^{\prime}$ denote this nonempty subset and call it a \emph{reduced block} of the quotient diagram of $(\textbf{a}, \textbf{b})$. If $\mathcal{C}_n ^{\prime}$ denotes the set of columns of $\mathcal{D}_n ^{\prime}$, one can then show that $\mathcal{C}_n ^{\prime}$  is the pairwise disjoint union of nonempty sets $\Omega_1,\dots,\Omega_r$ of adjacent essential columns which possess the following two properties:
\vspace{0.3cm}

$(i)$If $|\Omega_i| \geq 2$ and $C, C^{\prime}$ are adjacent essential columns of $\Omega_i$ then the exits a row of $\mathcal{D}(R_n)$ labeled by an element of $\Sigma_n$  which  passes through $C$ and $C^{\prime}$;

$(ii)$ If $i \not= j, C \in \Omega_i$, and $C^{\prime} \in \Omega_j$, then there does not exit a row of $\mathcal{D}(R_n)$ labeled by an element of $\Sigma_n$  which  passes through $C$ and $C^{\prime}$.
\vspace{0.3cm}

\noindent We will refer to each of the sets $\Omega_i, i \in [1, r]$, as a \emph{cell} of the reduced block $\mathcal{D}_n ^{\prime}$.

Let
\[
E_i= \bigcup_{C \in \Omega_i}\ \theta(C),\ i \in [1, r].
\]
It follows from (3.32) that whenever $\Sigma_n \not= \emptyset,$
\begin{equation*}
\Sigma_n=\bigcup_i E_i. \tag{3.33}
\] 
We note that the union in (3.33) is pairwise disjoint: if $s \in E_i \cap E_j$ for $i \not= j$ then there are columns $C \in \Omega_i$ and $C^{\prime} \in \Omega_j$ such that $s \in \Sigma_n \cap \theta(C) \cap \theta(C^{\prime})$, contrary to property $(ii)$ of the cells. 
\begin{lem}
\label{lem11}
If $\Sigma_n \not= \emptyset$ then $\{E_1,\dots, E_r\}$ is the set of equivalence classes of $\sim$ that are contained in $\Sigma_n$.  
\end{lem}

\emph{Proof}. Let $n\in [1, r]$. We first prove that $E_n$ is contained in an equivalence class of $\sim$. In order to do that, we note first that if $C$ is an essential column in $\mathcal{C}_n$ then it follows from Lemma 3.10 that $\theta(C)$ is contained in an equivalence class of $\sim$. We thus conclude from the linking property $(i)$ of the sections that $E_n$ is also contained in an equivalence class. Suppose next that $i \in E_n, j \in \Sigma \setminus \{i\}$, and $i \sim j$. As a consequence of the definition of $\sim$ and the proof of Lemma 3.9, there is a sequence $(s_1,\dots,s_m)$ of elements of $\Sigma_n$ connecting $i$ to $j$, with associated link sequence $(\theta(C_1),\dots, \theta(C_{m-1}))$ and $C_l \in \mathcal{C}_n^{\prime}$ for each $l \in [1, m-1]$. Now choose $C \in \Omega_n$ such that $i \in \theta(C)$. Then $i \in \theta(C) \cap \theta(C_1)$, hence it follows from property $(ii)$ of the cells that $C_1 \in \Omega_n$. An inductive argument which uses the fact that $\{s_l, s_{l+1}\} \subseteq \theta(C_l), l \in [1, m-1],$ now shows that $C_l \in \Omega_n$ for all $l \in [1, m-1]$. In particular, $j=s_m \in \theta(C_{m-1}) \subseteq E_n$, and so we conclude that $E_n$ is an equivalence class of $\sim$.

It follows that every element of  $\{E_1,\dots, E_r\}$ is an equivalence class of $\sim$ and because of (3.33), there are no others. $\    \Box$
\vspace{0.3cm}

Lemmas 3.9 and 3.11 now justify the following procedure for calculation of the equivalence classes of $\sim$: for any admissible $2k$-tuple,
\vspace{0.3cm}

(1) determine the reduced blocks of the quotient diagram, 

(2) locate the set of cells in each reduced block, and

(3) for each cell $\Omega$ in a reduced block, calculate the set formed by the labels of all rows of the associated block of the quotient diagram that pass through a column in $\Omega$.
\vspace{0.3cm}

\noindent  The parameter $\mu$ is then simply the total number of cells from all of the reduced blocks.

Next we investigate the effect of condition (3.9) on the cell structure in the reduced blocks. Let $\mathcal{D}_n^{\prime}$ be a reduced block; we list the cells $\Omega_1,\dots,\Omega_r$ of $\mathcal{D}_n^{\prime}$ as they are encountered as we move from left to right in  $\mathcal{D}(R_n)$ . Assume that $r \geq 2$ and let $\Omega_i$ and $\Omega_{i+1}$ be a pair of adjacent cells. A moment's reflection reveals the existence of a set $\{C_1,\dots, C_t\}$ of non-essential columns which possesses the following properties: 
\begin{equation*}
\theta(C_i) \cap \theta(C_{i+1}) \not= \emptyset\ \textrm{if}\ i\in [1, t-1], \tag{3.34}
\]
and if $C$ is the last column in $\Omega_i$ and $C^{\prime}$ is the first column in $\Omega_{i+1}$ as we proceed from left to right in $\mathcal{D}(R_n)$ then
\begin{equation*}
\theta(C) \cap \theta(C_1) \not= \emptyset \not= \theta(C^{\prime}) \cap \theta(C_t). \tag{3.35}
\]
Because $C_i$ is non-essential for all $i \in [1, t]$, we reason inductively via (3.34) to conclude that there is an integer $\sigma_0$ such that
\begin{equation*}
\sigma_i=\sigma_0,\ \textrm{for all}\ i \in \bigcup_{l=1}^t\ \theta(C_l). 
\tag{3.36}
\]
Now choose $i \in \theta(C) \cap \theta(C_1)$ and $j \in \theta(C^{\prime}) \cap \theta(C_t)$. As a consequence of (3.32), $i \in \Sigma_n$ and $j \in \Sigma_n$, hence it follows from property $(ii)$ of the cells that $i \not= j$. We also conclude from (3.35) and (3.36) that $\sigma_i=\sigma_0=\sigma_j.$ Thus, whenever $r \geq 2$, the map $i \rightarrow \sigma_i, i \in \Sigma$, is not injective, and so condition (3.9) cannot hold in this case. We conclude that if condition (3.9) is valid then there is exactly one cell, and hence also exactly one equivalence class of $\sim$ in $\Sigma _n$ whenever $\Sigma_n \not= \emptyset$. In particular, $\mu$ is the number of reduced blocks of the quotient diagram.

If we assume that all of the square-free parts of the coordinates of $\textbf{b}$ in $(\textbf{a}, \textbf{b})$ are distinct, then the calculation of the equivalence classes is even easier. In this case we have that $S(I)=I$ for all $I \in \Lambda(\mathcal{K})$ and so $ \Lambda^{\prime}(\mathcal{K})=\Lambda(\mathcal{K})$ and $\Sigma=\bigcup_i D_i$. 

We will now prove that $\{D_1,\dots,D_m\}$ is the set of equivalence classes of $\sim$. This will be true, by way of Lemma 3.9, if we show that for each $n \in [1, m]$ all elements of $D_n$ are $\sim$-equivalent. The following terminology will be of use in the argument : if $C$ is a column of a block $\mathcal{D}(R_n)$ in the quotient diagram of $(\textbf{a}, \textbf{b})$ then we will say that $C$ \emph{hangs from}  row $i$ \emph{in} $\mathcal{D}(R_n)$ if the top point of $C$ is in row $i$. We may assume with no loss of generality that for some integer $d$, the rows of $\mathcal{D}(R_n)$ are labeled from 1 to $d$, i.e., we take $D_n=[1, d]$. Observe now  that for each of the rows from 1 through $d-1$, there is at least one column hanging from that row which has cardinality at least 2; otherwise there is a coordinate of the gap sequence of $\mathcal{D}(R_n)$ that exceeds $s-1$, which is impossible by construction of the quotient diagram. Thus for all $i \in [1, d-1]$, there is a column $C$ of $\mathcal{D}(R_n)$ such that $[i, i+1] \subseteq \theta(C)$, hence $i \sim i+1$ for all $i \in [1, d-1]$.

It follows that whenever $(\textbf{a}, \textbf{b})$ is an admissible $2k$-tuple for which the square-free parts of the coordinates of $\textbf{b}$ are distinct and satisfy condition (3.9), the cardinality of $\bigcup_i D_i$, the number of blocks $m$ in the quotient diagram, the set $\mathcal{M}_1$, and the set $D_{n_0}$ in which is located the index $i_0$ that determines $\mathcal{M}_1$ as per (3.19) completely determine the density of $\Pi_+(\textbf{a}, \textbf{b})$ by means of the corresponding formulae of Theorem 3.7. In particular, the density of $\Pi_+(\textbf{a}, \textbf{b})$ is at most $2^{-m}$ whenever  $\mathcal{M}_1= \emptyset$ or $\mathcal{M}_1= \Phi(D_{n_0})$ and is at most $(2^m-1)/2^{2m-1}$, otherwise. This gives an interesting number-theoretic interpretation to the number of blocks in the quotient diagram. In fact, if $\mathcal{A}$ denotes the set of all admissible $2k$-tuples for which the square-free parts of the coordinates of $\textbf{b}$ are distinct and satisfy (3.9), $k \in [2, +\infty)$, and if  $m \in [1, +\infty)$ then Lemma 3.8 can be used to show that there exists infinitely many elements $(\textbf{a}, \textbf{b})$ of $\mathcal{A}$ such that the quotient diagram of  $(\textbf{a}, \textbf{b})$ has $m$ blocks and the density of $\Pi_+(\textbf{a}, \textbf{b})$ is $2^{-m}$ (respectively, $(2^m-1)/2^{2m-1})$. One can also show that if $\{l, n\} \subseteq [1, +\infty),$ with $l \geq 2n,$ then there are infinitely many elements $(\textbf{a}, \textbf{b})$ of $\mathcal{A}$ such that  the density of $\Pi_+(\textbf{a}, \textbf{b})$ is $2^{1-l}(2^n-1)$.

Another interesting class of examples, lying at an extreme opposite from the admissible $2k$-tuples, consists of the standard $2m$-tuples $(\textbf{a}, \textbf{b})$ which satisfy the following condition: if $B$ is the set of distinct values of the coordinates of $\textbf{b}$ and the sets $A(b), b \in B$, are defined as in the Introduction, then the sets $\{a/b: a \in A(b)\}, b \in B,$ are totally ordered by inclusion. In this case it follows that if $\mathcal{D}(R)$ is a block of the quotient diagram of $(\textbf{a}, \textbf{b})$ and $C$ is a column of $\mathcal{D}(R)$ then $K(C)=[\min K(C), |B|]$. The corresponding equivalence relation $\sim$ hence has only one equivalence class and $\Sigma=[1, |B|]$. It follows that whenever the square-free parts of the coordinates of $\textbf{b}$ satisfy (3.9), the density of $\Pi_+(\textbf{a}, \textbf{b})$ is $2^{1-|B|}$. 

 Let $(\mathbf{a}, \mathbf{b})$ be a fixed standard $2m$-tuple. We conclude our discussion in this paper by indicating how Theorem 3.7 needs to be modified so as to obtain the density of $\Pi_+(\mathbf{a}, \mathbf{b})$ when condition (3.9)  no longer holds, modulo the solution of a problem in enumerative combinatorics that we do not solve. In order to explain these modifications efficiently, we first recall that if $I \in \Lambda^{\prime}(\mathcal{K})$ then $\mathcal{S}(I)$ denotes the set $\{\pi(\sigma_i): i \in I\}$.

Next, let $\mathcal{M}_0$ and $\mathcal{M}_1$ be as defined by (3.19) and (3.20),  let $\Sigma$ be defined as before, and let $\mu$ denote the cardinality of $\Sigma/\sim$, also as before. We now define two combinatorial parameters $\alpha$ and $\beta$. First, consider the subset of $\mathcal{P}_{\emptyset}(\Lambda^{\prime}(\mathcal{K}), 2)$ consisting of all elements $\{P_1, P_2\}$ which satisfy the following condition: there exists a subset $U$ of
\[
\big\{\pi(\sigma_i): i \in \Sigma \big\} \cup \{\emptyset\}
\]
of odd cardinality such that the cardinality of
\[
U \cap \Big( \Big( \bigcup_{I \in P_1}\ \mathcal{S}(I) \Big) \cup \{\emptyset\} \Big)
\]
is even and the repeated symmetric difference of the elements of $U$ is empty. Let $\alpha$ denote the cardinality of this subset. Secondly, assume that $\mathcal{M}_1 \not= \emptyset$ and then consider the subset of  $\mathcal{P}_{\emptyset}(\mathcal{M}_0, 2)$ consisting of all elements  $\{P_1, P_2\}$ which satisfy the following condition: there exists a subset $U$ of 
\[
\big\{\pi(\sigma_i): i \in \Sigma \big\} 
\]
of odd cardinality such that the cardinality of
\[
U \cap \Big( \bigcup_{I \in \mathcal{M}_1 \cup P_1}\mathcal{S}(I) \Big) 
 \]
is even and the repeated symmetric difference of the elements of $U$ is empty. Let $\beta$ denote the cardinality of this subset.

Now, as before, let
\vspace{0.2cm}
\begin{equation*}
\varepsilon= \left\{\begin{array}{cc}1,\ \textrm{if $\mathcal{N}( \mathcal{M}_0, \mathcal{M}_1) \not= \emptyset$,}\\
0,\ \textrm{otherwise,}\end{array}\right.
\vspace{0.2cm}
\]
\[
n=\big| \bigcup_{i \in \Sigma}\ \pi(\sigma_i) \big|,
\vspace{0.2cm}
\]
\[
d=\ \textrm{the dimension of the linear span of $\{v(\pi(\sigma_i)): i \in \Sigma\} \setminus \{0\}$ in $F^n$.}
\vspace{0.2cm}
\]
 Then the density of $\Pi_+(\mathbf{a}, \mathbf{b})$  is
\[
2^{1-d}(2^{\mu-1}- \alpha),\  \textrm{if}\  \mathcal{M}_1= \emptyset.
\]

Assume next that $\mathcal{M}_1 \not= \emptyset$, let
\[
\Omega=\{\varpi \in \Sigma/\sim: \mathcal{M}_1 \cap \Phi(\varpi) \not= \emptyset \},
\]
and set
\[
\omega=|\Omega|.
\]
A straightforward modification of our previous reasoning then shows that 
\vspace{0.2cm}
\begin{equation*}
\big|\mathcal{P}_{\emptyset}(\mathcal{M}_0, 2) \big|= \left\{\begin{array}{cc}2^{\mu- \omega}-1,\ \textrm{if $ \mathcal{M}_1 \not= \bigcup_{\varpi \in \Omega} \Phi(\varpi)$,}\\
0,\ \textrm{if $ \mathcal{M}_1= \bigcup_{\varpi \in \Omega} \Phi(\varpi)$ and $\mu= \omega$,}\\
2^{\mu- \omega-1}-1, \ \textrm{if $ \mathcal{M}_1= \bigcup_{\varpi \in \Omega} \Phi(\varpi)$ and $\mu> \omega$,}\end{array}\right.
\vspace{0.2cm}
\]
and that $\varepsilon=1$ here if and only if $\mathcal{M}_1= \bigcup_{\varpi \in \Omega} \Phi(\varpi)$ and for all subsets $U$ of $\big\{\pi(\sigma_i): i \in \Sigma \big\} $ of odd cardinality such that the cardinality of
\[
U \cap \Big( \bigcup_{I \in \mathcal{M}_1}\mathcal{S}(I) \Big) 
 \]
is even, the repeated symmetric difference of the elements of $U$ is not empty. Hence in this case, the density of $\Pi_+(\mathbf{a}, \mathbf{b})$ is
\vspace{0.2cm}
\[
2^{-d}(2^{\mu-\omega+1}-2 \beta-1),\ \textrm{if $\emptyset \not= \mathcal{M}_1\not=\bigcup_{\varpi \in \Omega}  \Phi(\varpi)$}, 
\vspace{0.2cm}
\]
or
\vspace{0.2cm}
\[
2^{-d}(2^{\mu-\omega}-2 \beta+ \varepsilon-1),\ \textrm{if}\ \mathcal{M}_1=\bigcup_{\varpi \in \Omega} \Phi(\varpi).
\vspace{0.2cm}
\]

When $(\mathbf{a}, \mathbf{b})$ satisfies condition (3.9), $\alpha= \beta=0$ and $\omega=1$, hence Theorem 3.7 gives ``local" maximum values of the coefficient of $2^{-d}$ in the formulae for the density of $\Pi_+(\mathbf{a}, \mathbf{b})$, and also the location of these local maxima, as $(\mathbf{a}, \mathbf{b})$ varies throughout all standard $2m$-tuples, $m \in [2, +\infty)$. The calculation of the parameters $\alpha$ and $\beta$ is an interesting problem in enumerative combinatorics that we invite the curious reader to contemplate.

\section{Conclusions}
If $(\textbf{a}, \textbf{b})$ is a standard $2m$-tuple such that the square-free parts of the coordinates of $\textbf{b}$ satisfy condition (3.9), then we have calculated the asymptotic density of $\Pi_+(\textbf{a}, \textbf{b})$. This density is determined by the formulae stated in Theorem 3.7, and depends on the parameters given by  the cardinality of $\Sigma=\bigcup_{I \in \Lambda^{\prime}(\mathcal{K})}\  Z(I)$, the cardinality of $\Sigma/ \sim$, the set $\mathcal{M}_1=\{I \in \Lambda^{\prime}(\mathcal{K}): 1 \in  S(I) \}$, and the equivalence class in $\Sigma/ \sim $ which contains $i_0$, where $i_0$ is  the index  that determines $\mathcal{M}_1$ as per (3.19), whenever $\mathcal{M}_1$ is not empty. We also indicate how the density of $\Pi_+(\textbf{a}, \textbf{b})$ can be calculated when condition (3.9) is not satisfied.

\end{document}